# Symmetries in generalized Kaprekar's routine

**Fernando Nuez**

**Retired professor of Universidad Politécnica de Valencia (Spain)**

## Abstract

Algebraic relations are established that determine the invariance of the transformed number after several transformations. The restrictions that determine the group structure of these relationships are analyzed, as is the case of the Klein group. Parametric $K^r$ functions associated with the existence of cycles are presented, as well as the role of the number of their links in the grouping of numbers in higher order equivalence classes. For this we have developed a methodology based on binary equivalence relations and the complete parameterization of the Kaprekar routine using $K_i$ functions of parametric transformation.

## 1. Introduction

The process known as Kaprekar's routine consists in sorting the digits in a number in descending order. Then, the opposite arrangement takes place, i.e., in ascending order. These two arrangements are subtracted one from another, resulting in the image of the original number. In 1949 Datlatreya Kaprekar showed that when iterating such a process with 4-digit numbers, the result was always 6174.

This routine of transformation has always aroused the interest among fans of mathematical games (Gardner, 1975; Nishiyama, 2006). But additionally there is an important research among mathematicians specialized in Number Theory. As an example, some of these works are indicated in references. Much of this effort has been directed to determining constants and cycles in various numbering systems and number of digits, as well as some aspects of transformation trees such as maximum distances.

However, unsuspected relationships are being uncovered such as the connection with Mersenne primes (Yamagami, 2018) and there is a current of opinion that understands there are important gaps in the understanding of the Kaprekar process. We are interested in investigating the possible connections with Group Theory.

In this work we will approach the study of the transformation in base 10 with generality of the number of digits and the following questions will be addressed:

- What are the algebraic relations that must exist between the numbers for their images to coincide after one or several transformations? For example:
    - Why do 83246529 and 17487561 return the same transformed number?
    - Why do 8178382562 and 4774473809 return the same number after two transformations?

- Why do 5068069 and 3071934 return the same number after four transformations?

- What algebraic structures underpin these symmetries?
- What is the relation between cycles and the algebraic architecture of the transformation trees or graphs?

## 2. The general process

Let any number n of w digits $n = a_1 a_2 \ldots a_w$ belong to set $A_w \subset \mathbb{Z}$, excluding numbers with identical digits. Let numbers of less digits be included, as long as they are completed by adding zeros to their left. Let us use the decimal numeral system, so that $0 \leq a_i \leq 9$, i= 1,...w

Let $O_d$ be an operator that sorts the digits of a number in descending order: $x = O_d (n) = x_1 \ x_2 \ldots x_w$, $x_i \geq x_{i+1}$ and let $O_u$ be another operator which does the opposite: $y = O_u (n) = x_w \ x_{w-i} \ldots x_1$

Generalized Kaprekar's routine is formalized through an operator K such that

$$K (n) = n' = x-y = O_d (n) - O_u (n) \qquad [1]$$

For example, if n = 83246529, $O_d$ (n) = x = 98654322

$O_u$ (n) = y = 22345689, n' = K (n) = x-y = 76308633

The iteration of the process r times is represented by

$K^r (n) = K^{r-1}[K(n)] \quad r \geq 1$

agreeing that $K^0$ (n) = n and writing $K^1 = K$ for convenience.

In our example $K^2$ (83246529) = 84326652

Arranging the digits through $O_d$ allows to define the following parameters, resulting from subtracting the digits in x that are symmetrical compared to the central value.

$\alpha^s = x_s - x_{w-s+1} \quad 1 \leq s \leq h$, h = w/2 if $w = \dot{2}$, h = (w-1)/2, if $w = \dot{2}+1$; $0 < \alpha^1 \leq 9$, $0 \leq \alpha^s \leq 9$, $2 \leq s \leq h$ [2]

$\alpha^s \geq \alpha^{s+1}$, s = 1,...h-1 [3]

must be verified as a consequence of the arrangements [1] forced by the routine.

Thus, for w = 4, $\alpha^1 = x_1\text{-}x_4$, $\alpha^2 = x_2\text{-}x_3$. For simplicity, these parameters will be called $\alpha = \alpha^1$ and $ß = \alpha^2$ in further examples. In a similar way, the same parameters exist in $A_5$.

As will be shown below, parameters $\alpha^s$ uniquely determine the image of a number. Consequently, any number is characterized by their parameters, which are represented as follows:

$p (n) = (\alpha^1 \ \alpha^2 \ldots \alpha^h) = \boldsymbol{\alpha}$, $n \in A_w$

When the situation is unambiguous, such parameters will be written together, as a number of h digits.

The transformation functions are the following:

**1) If $0 < \alpha^s \leq 9,\ s = 1, 2, \ldots h$**

**1.1) $w = \dot{2},\quad h = w/2$**

$f_I(\alpha^1\ \alpha^2 \ldots \alpha^h) = (x_1\text{-}x_w\quad x_2\text{-}x_{w-1} \ldots x_{w-1}\text{-}x_2\quad x_w\text{-}x_1) =$

$= (\alpha^1\ \alpha^2 \ldots \alpha^h\quad -\alpha^h \ldots -\alpha^2\ \ -\alpha^1) =$

$= (\alpha^1\ \alpha^2 \ldots \alpha^{h-1}\quad \alpha^h\text{-}1\quad 9\text{-}\alpha^h\quad 9\text{-}\alpha^{h-1}\ 9\text{-}\alpha^2\ \ 10\text{-}\alpha^1)$ [4a]

So for 6-digit numbers (w = 6)

$f_I(\alpha^1\ \alpha^2\ \alpha^3) = (\alpha^1\ \alpha^2\ \alpha^3 - 1\quad 9\text{-}\alpha^3\ \ 9\text{-}\alpha^2\ \ 10\text{-}\alpha^1)$

which can be written as follows to avoid using superscripts:

$f_I(\alpha\ ß\ \gamma) = (\alpha\ \ ß\ \ \gamma\text{-}1\ \ 9\text{-}\gamma\ \ 9\text{-}ß\ \ 10\text{-}\alpha),\ \alpha = \alpha^1 > 0,\ ß = \alpha^2 > 0\ ,\ \gamma = \alpha^3 > 0$

For instance, what is the image of 631764?

$O_d\,(631764) = 766431,\ \alpha = 7\text{-}1 = 6,\ ß = 6\text{-}3 = 3, \gamma = 6\text{-}4 = 2$

$p\,(631764) = 632,\quad f_I\,(632) = 631764$

K (631764) = 631764

Apparently, this number transforms into itself. It is a transformation or balance constant in $A_6$.

However, note that any number with these same parameters with yield the same image $n_E$. For example,

$m = 632000,\quad p\,(m) = 632,\quad m' = K\,(m) = n_E$

Any permutation P of the digits in

$(6+a\ \ 3+b\ \ 2+c\ \ a\ \ b\ \ c),\ \ 0 \leq a \leq 3,\ \ 0 \leq b \leq 6,\ \ 0 \leq c \leq 7,\ \ a \geq b \geq c,\ \ a > c$

will become $n_E$, $K\,[P\,(6+a\ \ 3+b\ \ 2+c\ \ a\ \ b\ \ c)] = 631764 = n_E$

Digits $a'_s$ in those numbers transformed by $f_I$ must satisfy the following restrictions:

$\sum_{s=1}^{w} a'_s = \dot{9} \rightarrow K(n) = \dot{9};\ a'_1 + a'_w = 10;\ a'_s + a'_{w-s+1} = 9,\ \ 1 < s < h$

$a'_h + a'_{h+1} = 8$ [4a.1]

In the case of $A_6$ : $a'_1 + a'_6 = 10,\ a'_2 + a'_5 = 9,\ a'_3 + a'_4 = 8,\ \sum_{s=1}^{6} a_s = 27$

**1.2) $w = \dot{2}+1,\quad h = (w\text{-}1)/2$**

In this case, function $f_I$ is similar to the last one, but it includes a 9 as a central digit in the transformed number

$f_I(\alpha^1\ \alpha^2 \ldots \alpha^h) = (\alpha^1\ \ \alpha^2 \ldots \alpha^{h-1}\ \ \alpha^h\text{-}1\quad 9\quad 9\text{-}\ \alpha^h\ \ 9\text{-}\ \alpha^{h-1} \ldots 9\text{-}\alpha^2\quad 10\text{-}\ \alpha^1)$ [4b]

For instance, the image of a telephone number

n = 34326714825 would be

$O_d = 87654433221,\ \alpha^1 = 7,\ \ \alpha^2 = 5,\ \ \alpha^3 = 4,\ \ \alpha^4 = 2,\ \ \alpha^5 = 1$

$f_I\,(75421) = (7\ \ 5\ \ 4\ \ 2\ \ 1\text{-}1\ \ 9\ \ 9\text{-}1\ \ 9\text{-}2\ \ 9\text{-}4\ \ 9\text{-}5\ \ 10\text{-}7) = 75420987543$

Notably, for 7-digit numbers

$f_I(\alpha\ ß\ \gamma) = (\alpha\ ß\ \ \gamma\text{-}1\ \ 9\ \ 9\text{-}\ \gamma\ \ 9\text{-}ß\ \ 10\text{-}\alpha)$

The digits in the transformed number must satisfy

$\sum_{s=1}^{w} a'_s = \dot{9} \rightarrow K(n) = \dot{9};\ a'_1 + a'_w = 10;\ a'_s + a'_{w-s+1} = 9,\ \ 1 < s < h$

$a_h + a_{w-h+1} = 8,\ \ a_{h-1} = 9$ [4b.1]

**2) If $0 < \alpha^s \leq 9$, $1 \leq s < r$ and $\alpha^s = 0$, $s \geq r$ (r>2)**

**2.1 ) $w = \dot{2}$, $h = w/2$**

Basic function $f_2$ is similar to [4a] but adding $v = 2(h+1-r)$ 9 digits to the middle part (two nines per each null parameter)

$$f_2(\alpha^1 \ldots \alpha^{r-1}\ 0\ \ 0 \ldots 0) = (\alpha^1 \ldots \alpha^{r-2}\ \ \alpha^{r-1}-1\ \ \overbrace{9 \ldots 9}^{v}\ \ 9-\alpha^{r-1}\ \ 9-\alpha^{r-2} \ldots 9-\alpha^2\ \ 10-\alpha^1) \quad [5a]$$

thus satisfying the conditions

$\sum_{s=1}^{w} a'_s = \dot{9} \rightarrow K(n) = \dot{9}$; $a'_1 + a'_w = 10$; $a'_s + a'_{w-s+1} = 9$, $2 \leq s \leq r-2$

$a'_{r-1} + a'_{w-r+2} = 8$, $a'_r = a'_{r+1} = \ldots = a'_{w-s+1} = 9$ [5a.1]

Specifically, for $A_6$

$f_2(\alpha\ \beta\ 0) = (\alpha\ \ \beta-1\ \ 9\ \ 9\ \ 9-\beta\ \ 10-\alpha)$

$\sum a'_s = 36 \rightarrow K(n) = \dot{9}$; $a'_1 + a'_6 = 10$; $a'_2 + a'_5 = 8$, $a'_3 = a'_4 = 9$

For example, let n = 549945 have the image

$O_d$ (n) = 995544, p (n) = 550, $f_2$ (550) = 549945, the original number itself. This is another constant in $A_6$.

**2.2 ) $w = \dot{2}+1$, $h = (w-1)/2$**

Function $f_2$ is similar to the last one but adds an extra 9 as a middle digit in the transformed number

$$f_2(\alpha^1 \ldots \alpha^{r-1}\ 0\ \ 0 \ldots 0) = (\alpha^1 \ldots \alpha^{r-2}\ \ \alpha^{r-1}-1\ \ \overbrace{9 \ldots 9}^{v+1}\ \ 9-\alpha^{r-1}\ \ 9-\alpha^{r-2} \ldots 9-\alpha^2\ \ 10-\alpha^1) \quad [5b]$$

adding to conditions [5a.1], $a'_{h+1} = 9$ [5b.1]

For $A_7$ and $\gamma = 0$, the result is

$f_2(\alpha\ \beta\ 0) = (\alpha\ \ \beta-1\ \ 9\ \ 9\ \ 9\ \ 9-\beta\ \ 10-\alpha)$ with

$a'_1 + a'_7 = 10$, $a'_2 + a'_6 = 8$, $a'_3 = a'_4 = a'_5 = 9$

**3) If $0 < \alpha^1 \leq 9$ and $\alpha^s = 0$, $s \geq 2$**

$$f_3(\alpha^1\ \ \overbrace{0 \ldots 0}^{h-1}) = (\alpha^1-1\ \ \overbrace{9 \ldots 9}^{w-2}\ \ 10-\alpha^1) \quad [6]$$

which implies the following restrictions for the transformed numbers

$a'_1 + a'_w = 9$, $a'_2 = \ldots = a'_{w-1} = 9$ [6.1]

These two expressions are valid both for even and odd number of digits

For $A_6$ the result is

$f_3(\alpha\ 0\ 0) = (\alpha-1\ \ 9\ \ 9\ \ 9\ \ 9\ \ 10-\alpha)$

And for $A_7$

$f_3(\alpha\ 0\ 0) = (\alpha-1\ \ 9\ \ 9\ \ 9\ \ 9\ \ 9\ \ 10-\alpha)$

Functions $f_1$, $f_2$ and $f_3$ are similar to the ones developed by Prichett et al. (1981). These functions are an epimorphism of $A_W$ in $B_W \subset A_W$ which comprises those multiples of 9 that verify some requirement [5a.1], [5b.1] or [6.1] and that have antiimage.

## 3 Parametric functions

The previous development returns

$p(m) = p(n) \leftrightarrow K(m) = K(n)$ [7]

which suggests analyzing the transformation process in parametric terms. To this end, it is necessary to know the parameters of the transformed numbers.

One of the questions in the introduction has thus been answered. Both numbers 83246529 and 17487561 return the same number after their transformations because both have the same parameters (7631).

Let us specify the notation used so far

- K, in a broad sense, to refer to the transformation. Its argument is a number and its image is the transformed number K (n) = n'
- $f_i$ in a specific sense. These act upon parameters and yield the numerical image of all numbers whose parameters are those specified.
  $p(n) = \boldsymbol{\alpha}$, $f_i(\boldsymbol{\alpha}) = n'$, $\boldsymbol{\alpha} = \alpha^1\ \alpha^2 \ldots \alpha^h$
- K can also be used broadly, by extension, to mean the parameters of some transformed number, thus transforming parameters into parameters
  $K(\boldsymbol{\alpha}) = \boldsymbol{\alpha}'$, $\boldsymbol{\alpha}' = p(n')$

This is the biggest problem in transformations. The digits in the transformed number are not necessarily arranged, so it is not possible to determine its parameters right away. Thus, operators $O_d$ and $O_u$, which establish a permutation $P_i$ that sorts the digits, must operate. $O_d(n') = P_i(n')$. If $P_i$ is identified, then $\boldsymbol{\alpha}'$ is directly discovered. A specific function $K_i(\boldsymbol{\alpha}) = \boldsymbol{\alpha}'$ can therefore be associated with each permutation $P_i$.

This is a tedious process, since the possible permutations will grow factorially as the number of digits increases:

- With w =2 and w = 3 there are just 2 functions $K_i(\alpha)$
- With w = 4 and w = 5 there are 11 functions $K_i(\alpha\ ß)$ and 2 functions $K_i(\alpha\ 0)$
- With w = 6 and w = 7 there are 2 functions $K_i(\alpha\ 0\ 0)$, 11 functions $K_i(\alpha\ ß\ 0)$ and 117 functions $K_i(\alpha\ ß\ \gamma)$
- With w and w+1 digits, the number of functions $K_i(\alpha^1\ \alpha^2 \ldots \alpha^s \ldots \alpha^h)$, $\alpha_s > 0$, s =1, …h approaches w!/h! by default. This is because in $O_d(n')$ a permutation with $\alpha^{r-s}$ to the left of $\alpha^r$ is, under [3], only possible if $\alpha^r = \alpha^{r-s}$.

  Let us see some examples of these functions.

1) **Functions $K_i$ based on [4a]**
   **$w = \dot{2}$, h = w/2, $f_1(\alpha)$ defined in [4a] and $P_1$ (1 2 3 …w)**
   $P_1$ is defined as the digit arrangement described in [4a]
   $K_1(\boldsymbol{\alpha}) = [\alpha^1 - (10-\alpha^1)\ \ \alpha^2 - (9-\alpha^2) \ldots \alpha^{h-1} - (9-\alpha^{h-1})\ \alpha^h - 1 - (9-\alpha^h)] =$
   $= [2\alpha^1 - 10\ \ 2\alpha^2 - 9 \ldots 2\alpha^{h-1} - 9\ 2\alpha^h - 10]$, [8a]

   The domain of existence of this function $K_1(\boldsymbol{\alpha}) = \boldsymbol{\alpha}'$ is determined by the restrictions imposed by $O_d(n')$, the digits being sorted in descending order [1] in permutation $P_1$ $\alpha^1 \geq \alpha^2, \ldots, \alpha^h-1 \geq 9-\alpha^h, 9-\alpha^h \geq 9-\alpha^{h-1}, \ldots, 9-\alpha^3 \geq 9-\alpha^2, 9-\alpha^2 \geq 10-\alpha^1$
which implies $\alpha^h \geq 5$, and under [3], $\alpha^s \geq 5$, $s \geq 1$, y $\alpha^1 \geq \alpha^2+1$
Briefly, the existence conditions for $K_1$ (**α**) are
$w = \dot{2}; 5 < \alpha^1 \geq \alpha^2+1; 5 \leq \alpha^s \leq 9, \; s \geq 2$ [8a.1]
For instance,
n = 181771978221, w = 12, h = 6, p(n) = 877655 = **α** which satisfies [8a.1
$K_1$ (877655) = 655310 and indeed, $f_1$(n)=877654443222=n', p(n')= 655310=**α**'
Note that this **α**' does not satisfy [8a.1], so it cannot be transformed through $K_1$. We need a function $K_i$ based on [5a], with $\alpha^h = \alpha^6 = 0$ (see [15] in the last case in section 3).

**w = $\dot{2}$, h = w/2, $f_1$ (α) defined in [4a] and $P_2$ (h+1 …w 1 2 … h)**
$O_d$ (n') = $( 9-\alpha^h \;\; 9-\alpha^{h-1} \ldots 9-\alpha^2 \;\; 10-\alpha^1 \;\; \alpha^1 \; \alpha^2 \ldots \alpha^{h-1} \; \alpha^h-1)$
$K_2$ (**α**) = $(10-2\alpha^h \;\; 9-2\alpha^{h-1} \ldots 9-2\alpha^2 \;\; 10-2\alpha^1)$
$w = \dot{2}; 0 < \alpha^1 \leq 5; \alpha^1 \geq \alpha^2+1; 0 < \alpha^s \leq 4, s=2 \ldots h, 2\alpha^h \leq 2\alpha^{h-1}+1$ [9a]
For instance, in $A_6$
$K_2$ (541) = 810

**w = $\dot{2}$, h = w/2 and $P_3$ (h+1 1 h+2…w 2 3… h) in [4a]**
$O_d$ (n') = $( 9-\alpha^h \;\; \alpha^1 \; 9-\alpha^{h-1} \;\; 9-\alpha^{h-2} \ldots 9-\alpha^3 \;\; 9-\alpha^2 \;\; 10-\alpha^1 \;\; \alpha^2 \ldots \alpha^{h-1} \; \alpha^h-1)$
$K_3$ (**α**) = $(10-2\alpha^h \;\; \alpha^1-\alpha^{h-1} \;\; 9-\alpha^{h-1}-\alpha^{h-2} \ldots 9-\alpha^3-\alpha^2 \;\; \alpha^1-\alpha^2-1)$
$w=\dot{2}$, h=w/2; $1 \leq \alpha^1+\alpha^2 \leq 10$; $\alpha^s+\alpha^{s+1} \leq 9$, s= 2,3,…h-2; $\alpha^1 \geq \alpha^2+1$; $\alpha^1+\alpha^{h-1} \geq 9$; $\alpha^1+\alpha^h \leq 9$; $\alpha^h \leq 10+\alpha^{h-1}-2\alpha^h$; $1 \leq \alpha^h \leq 4$ [10]
This is an important function that generates a family of transformation constants.

## 2) <u>Functions $K_i$ based on [4b]</u>

The permutations refer to the digits to the right of 9 in [4b]

**w = $\dot{2}$+1, h = (w-1)/2, $f_1$ (α) defined in [4b] and $P_1$**
$O_d$ (n') = $( 9 \;\; \alpha^1 \;\; \alpha^2 \ldots \alpha^{h-1} \;\; \alpha^h-1 \;\; 9-\alpha^h \;\; 9-\alpha^{h-1} \ldots 9-\alpha^2 \; 10-\alpha^1)$
The presence of a 9 completely changes function $K_1$
$K_1$ (**α**) = $(\alpha^1-1 \;\; \alpha^1+\alpha^2-9 \;\; \alpha^2+\alpha^3-9 \ldots \alpha^{h-1}+\alpha^h-9)$ [8b]
with the existence conditions
$w = \dot{2}+1, 5 < \alpha^1 \geq \alpha^2+1; \;\; 5 \leq \alpha^s \leq 8, s > 1; \alpha^s+\alpha^{s+1} \geq 9, s = 1, \ldots h-1$ [8b.1]
For example, n = 8650000, p (n) = 865 = **α**
n' = $f_1$ (n) = 8649432, $O_d$ (n') = 9864432, p (n') = 752 = **α**'
**α** = 865 satisfies the conditions [8b.1] and, in fact,
$K_1$ (**α**) = (8-1 8+6-9 6+5-9) = 752

**w = $\dot{2}$+1, h = (w-1)/2, and $P_2$ in [4b]**
$O_d$ (n') = $( 9 \;\; 9-\alpha^h \;\; 9-\alpha^{h-1} \ldots 9-\alpha^2 \;\; 10-\alpha^1 \;\; \alpha^1 \;\; \alpha^2 \cdots \alpha^{h-1} \;\; \alpha^h-1)$
$K_2$ (**α**) = $(10-\alpha^h \;\; 9-\alpha^h-\alpha^{h-1} \ldots 9-\alpha^3-\alpha^2 \;\; 9-\alpha^2-\alpha^1)$ [9b]
$w = \dot{2}+1$, h=(w-1)/2; $0 < \alpha^1 \leq 5$; $\alpha^1 \geq \alpha^2+1$; $\alpha^{s+1}+\alpha^s \leq 9$, s=1, 2, ... h-1
$\alpha^{h-1} \geq 1$, $1 < \alpha^h < 5$
For example, n = 2515324 p(n) = 432
w = 7, h = 3 por [4b], n' = $f_1$ (n) = 4319766, p (n') = 842
and in fact, since 432 satisfies the existence conditions of $K_2$
$K_2$ (432) = ( 10-2 9-2-3 9-3-4) = 842

**$w = \dot{2}+1$, $h = (w-1)/2 \geq 4$, $P_4$ (1 h+2 2 h+3 … h+s 3 4 … s 2h 2h+1 h), s = h-1 en [4b]**

$O_d$ (n') = $(9\ \ \alpha^1\ \ 9\text{-}\alpha^h\ \ \alpha^2\ \ 9\text{-}\alpha^{h-1} \ldots 9\text{-}\alpha^3\ \ \alpha^3 \ldots \alpha^{h-1}\ \ 9\text{-}\alpha^2\ \ 10\text{-}\alpha^1\ \ \alpha^h\text{-}1)$

$K_4(\boldsymbol{\alpha}) = (10\text{-}\alpha^h\ \ 2\alpha^1\text{-}10\ \ \alpha^2\text{-}\alpha^h\ \ \alpha^2\text{-}\alpha^{h-1}\ \ 9\text{-}\alpha^{h-2}\text{-}\alpha^{h-1} \ldots 9\text{-}\alpha^4\text{-}\alpha^3)$ [10b]

$h \geq 4;\ 9 \leq \alpha^1 + \alpha^h \leq 11,\ \alpha^2+\alpha^h \leq 9,\ \alpha^2+\alpha^{h-1} \geq 9,\ \alpha^3 \leq 4$

Por ejemplo, en w = 15, h = 7

$f_1(\boldsymbol{\alpha}) = (\alpha^1\ \alpha^2 \ldots \alpha^6\ \ \alpha^7\text{-}1\ \ 9\ \ 9\text{-}\alpha^7\ \ 9\text{-}\alpha^6 \ldots 9\text{-}\alpha^2\ \ 10\text{-}\alpha^1) = n'$

$O_d$ (n') = $(9\ \ \alpha^1\ \ 9\text{-}\alpha^7\ \ \alpha^2\ \ 9\text{-}\alpha^6 \ldots 9\text{-}\alpha^3\ \ \alpha^3 \ldots \alpha^6\ \ 9\text{-}\alpha^2\ \ 10\text{-}\alpha^1\ \ \alpha^7\text{-}1)$

$K_4(\boldsymbol{\alpha}) = (10\text{-}\alpha^7\ \ 2\alpha^1\text{-}10\ \ \alpha^2\text{-}\alpha^7\ \ \alpha^2\text{-}\alpha^6\ \ 9\text{-}\alpha^5\text{-}\alpha^6\ \ 9\text{-}\alpha^4\text{-}\alpha^5\ \ 9\text{-}\alpha^3\text{-}\alpha^4)$

$K_4$ (8643332) = 864332 therefore this number is a parametric constant in w = 15, with the numeric constant of 15-digits, 864333197666532.

Similarly, for 9-digits

K4 (α ß γ δ) = (10-δ 2α-10 ß-δ ß- γ), $K_4$ (8642) = 8642, which is a parametric constant in $A_8$

$K_4(\boldsymbol{\alpha})$, defined in a general way at $w = \dot{2}+1$, is an important parametric function from which 1-cycles derive.

### 3) <u>Functions $K_i$ based on [5a]</u>

Functions $K_i$ based on [5a], h = w/2, are complex, as they depend on the number of null parameters compared to the number of non-null parameters. If $r \leq 1+h/2$, $\alpha^r = 0$, then the functions' structure remains, since [5] provides as many 9 as it provides varying digits. If there are more non-null than null parameters, r > 1+h/2, then functions $K_i$ will vary with each r increment.

Let us see some examples,

- **$w = \dot{2}$, $r \leq 1+h/2$, $\alpha^3 = 0$, $h \geq 4$ in [5a] and $P_1$**

$n' = f_2(\alpha^1\ \alpha^2\ 0\ \overbrace{\ldots\ldots}^{h-2}\ 0) = (\alpha^1\ \alpha^2\text{-}1\ \ 9\ \overbrace{\ldots\ldots}^{w-4}\ 9\ \ 9\text{-}\alpha^2\ \ 10\text{-}\alpha^1)$

If using permutation $P_1$ (1234) to the right of numbers 9 in $O_d$ (n')

$O_d(n') = (9\ \overbrace{\ldots\ldots}^{w-4}\ 9\ \ \alpha^1\ \ \alpha^2\text{-}1\ \ 9\text{-}\alpha^2\ \ 10\text{-}\alpha^1)$

the associated function is

$K_{21}(\alpha^1\ \alpha^2\ 0\ \overbrace{\ldots\ldots}^{h-2}\ 0) = (\alpha^1\text{-}1\ \ \alpha^2\ \ 10\text{-}\alpha^2\ \ 9\text{-}\alpha^1\ \ 0\ \overbrace{\ldots\ldots}^{h-4}\ 0),\ 6 \leq \alpha^1 \leq 9,$

$5 \leq \alpha^2 \leq 9,\ \alpha^1 \geq \alpha^2+1$ [11a]

For example: $K_{21}$ (85000) = (75510)

- **$w = \dot{2}$, $r \geq 1+h/2$, $\alpha^3 = 0$, $h = 3$ in [5a] and $P_1$**

$n' = f_2(\alpha^1\ \alpha^2\ 0) = (\alpha^1\ \alpha^2\text{-}1\ \ 9\ \ 9\ \ 9\text{-}\alpha^2\ \ 10\text{-}\alpha^1)$

$O_d(n') = (9\ \ 9\ \ \alpha^1\ \ \alpha^2\text{-}1\ \ 9\text{-}\alpha^2\ \ 10\text{-}\alpha^1)$

$K_{22}(\alpha^1\ \alpha^2\ 0) = (\alpha^1\text{-}1\ \ \alpha^2\ \ \alpha^1\text{-}\alpha^2\text{+}1),\ 6 \leq \alpha^1 \leq 9,\ 5 \leq \alpha^2 \leq 9,\ \alpha^1 \geq \alpha^2+1,\ 2\alpha^2 \geq \alpha^1+1$ [12]

a function which differs from [11a]

For example: $K_{22}$ (850) = 754

- **$w = \dot{2}$, $r \leq 1+h/2$, $\alpha^3 = 0$, $h \geq 4$ in [5a], $P_2$ (1423)**

$O_d\,(n') = (\overbrace{9\ \ldots\ldots\ 9}^{w-4}\ \ \alpha^1\ \ 10-\alpha^1\ \ \alpha^2-1\ \ 9-\alpha^2)$

$K_{23}(\alpha^1\ \alpha^2\ \overbrace{0\ \ldots\ldots\ 0}^{h-2}) = (\alpha^2\ \ 10-\alpha^2\ \ \alpha^1-1\ \ 9-\alpha^1\ \ \overbrace{0\ \ldots\ldots\ 0}^{h-4})$, $5 \le \alpha^1 \le 9$, $5 \le \alpha^2 \le 9$, $\alpha^1+\alpha^2 \le 11$ [13]

For example: $K_{23}$ (5500) = 5544

- **$w = \dot{2}$, $r > 1+h/2$, $\alpha^3 = 0$, $h = 3$ in [5a], $P_2$**
  $O_d\,(n') = (9\ \ 9\ \ \alpha^1\ \ 10-\alpha^1\ \ \alpha^2-1\ \ 9-\alpha^2)$
  $K_{24}\,(\alpha^1\ \alpha^2\ 0) = (\alpha^2\ \ 10-\alpha^2\ \ 2\alpha^1-10)$, $5 \le \alpha^1 \le 9$, $5 \le \alpha^2 \le 9$, $\alpha^1+\alpha^2 \le 11$ [14]
  which operates only on parameters (550) and (650).
  With $r \le 1+h/2$, functions differ from those with $r > 1+h/2$

- **$w = \dot{2}$, $\alpha^h = 0$, $h = w/2 = 6$ in [5a], $P_3$ = (6 1 7 2 3 8 9 10 4 5)**
  As an example of applying [8a] we found
  n' = 877654443222, **α**' = p (n') = 655310. Such **α**' cannot be transformed with $K_1$. It is a case w = 12, h = 6, $\alpha^6 = 0$ which can be transformed with the function
  $K_{25}\,(\alpha\ ß\ \gamma\ \delta\ \varepsilon\ 0) = (10-\varepsilon\ \ 9-\delta\ \ \alpha-\varepsilon-1\ \ \alpha+ß-9\ \ \gamma-\delta\ \ ß-\gamma)$,
  $\alpha > 1$, $\alpha \ge ß+1$, $\alpha+ß \ge 9$, $9 \le \alpha+\delta \le 10$, $\alpha+\varepsilon \le 9$
  $\alpha+\delta-\varepsilon \le 10$, $ß+\delta \le 9$, $ß+\varepsilon \le 8$, $2\gamma \ge ß+\delta$, $\gamma \ge 5$, $\delta \ge \varepsilon+1$, $\alpha+ß+\delta-\gamma \ge 9$, $\varepsilon \ge 1$ [15]
  This function is associated with the aforementioned permutation $P_3$, in a specific case of [5a]
  $f_2\,(\alpha\ ß\ \gamma\ \delta\ \varepsilon\ 0) = (\alpha\ ß\ \gamma\ \delta\ \varepsilon-1\ \ 9\ \ 9\ \ 9-\varepsilon\ \ 9-\delta\ \ 9-\gamma\ \ 9-ß\ \ 10-\alpha)$,
  $f_2$ (655310) = 655309986444, $K_{25}$ (655310) = 964220

## 4) Functions $K_i$ based on [5b]

As previously stated, when the amount of digits is odd, $w = \dot{2}+1$, $h = (w-1)/2$, function $f_2$ is similar to the one corresponding for numbers with $w = \dot{2}$, but it includes one more 9 as a middle digit in the transformed number.

Such addition does not alter functions $K_i$ if $r \le 1+h/2$, $\alpha^r = 0$, although it does change the last parameters of **α**' if $r > 1+h/2$.

In parallel with the cases considered in 3)

- **$\alpha^3 = 0$, w = 7, h = 3 in [5b] and $P_1$ (1234)**
  $K_{26}\,(\alpha^1\ \alpha^2\ 0) = (\alpha^1-1\ \ \alpha^2\ \ 10-\alpha^2)$, $5 \le \alpha^1 \le 9$, $5 \le \alpha^2 \le 9$, $\alpha^1 \ge \alpha^2+1$ [16]
  a closer function to $K_{21}$ than to $K_{22}$
- **$\alpha^3 = 0$, w = 7, h = 3 in [5b] and $P_2$ (1423)**
  $K_{27}\,(\alpha^1\ \alpha^2\ 0) = (\alpha^2\ \ 10-\alpha^2\ \ \alpha^1-1)$, $5 \le \alpha^1 \le 9$, $5 \le \alpha^2 \le 9$, $\alpha^1+\alpha^2 \le 11$ [17]
  also a closer expression to $K_{23}$ than to $K_{24}$

## 5) Functions $K_i$ based on [6]

There are only two functions Ki based on [6]

$K_{31}\,(\alpha\ \ \overbrace{0\ \ldots\ldots\ 0}^{h-1}) = (\alpha-1\ \ 10-\alpha\ \ \overbrace{0\ \ldots\ldots\ 0}^{h-2})$, $6 \le \alpha \le 9$ [18]

$$K_{32}(\alpha\ 0\ \underbrace{\ldots\ldots}_{h-1}\ 0) = (10-\alpha\ \ \alpha-1\ \ 0\ \underbrace{\ldots\ldots}_{h-2}\ 0),\ 0 \le \alpha \le 5 \qquad [19]$$

valid both for $w = \dot{2}$ with $h = w/2$ and for $w = \dot{2}+1$ with $h = (w-1)/2$

Tables 1 and 2 show some functions $K_i$ of the different types in $A_6$ ($w = 6$) and $A_7$ ($w = 7$).

Table 1. Some functions $K_i$ in $A_6$, $w = 6$, $h = 3$

| Permutation O (n') | Associated $K_i$ ($\boldsymbol{\alpha}$) | Domain of existence |
|---|---|---|
| $P_1$ (123456) | $K_1(\alpha\ ß\ \gamma) = (2\alpha-10\ \ 2ß-9\ \ 2\gamma-10)$ | $6 \le \alpha \le 9$, $5 \le ß \le 9$, $5 \le \gamma \le 9$, $\alpha \ge ß+1$ |
| $P_2$ (456123) | $K_2(\alpha\ ß\ \gamma) = (10-2\gamma\ 9-2ß\ \ 10-2\alpha)$ | $0 < \alpha \le 5$, $0 < ß \le 4$, $0 < \gamma \le 4$, $\alpha \ge ß+1$, $2\gamma \le 2ß+1$ |
| $P_3$ (415623) | $K_3(\alpha\ ß\ \gamma) = (10-2\gamma\ \alpha-ß\ \ \alpha-ß-1)$ | $9 \le \alpha+ß \le 10$, $\alpha \ge ß+1$, $1 \le \gamma \le 4$<br>$\alpha+\gamma \le 9$, $\alpha \le 10+ß-2\gamma$ |
| $P_4$ (142563) | $K_4(\alpha\ ß\ \gamma) = (\alpha-\gamma+1\ \ \alpha-\gamma-1\ \ 2ß-9)$ | $5 \le \alpha \le 9$, $5 \le ß \le 8$, $1 \le \gamma \le 6$<br>$\alpha \ge ß+1$, $9 \le \alpha+\gamma \le 11$, $\alpha \ge 2ß+\gamma-8$ |
| $P_5$ (145263) | $K_5(\alpha\ ß\ \gamma) = (\alpha-\gamma+1\ \ \alpha-\gamma-1\ \ 9-2ß)$ | $6 \le \alpha \le 9$, $1 \le ß \le 4$, $\alpha+ß \ge 10$<br>$9 \le \alpha+\gamma \le 11$, $\alpha \ge 10-2ß+\gamma$ |
| $P_6$ (412563) | $K_6(\alpha\ ß\ \gamma) = (10-2\gamma\ 2\alpha-10\ \ 2ß-9)$ | $6 \le \alpha \le 9$, $5 \le ß \le 9$, $1 \le \gamma \le 3$<br>$\alpha \ge ß+1$, $\alpha+\gamma \le 9$ |
| $P_7$ (124536) | $K_7(\alpha\ ß\ \gamma) = (2\alpha-10\ \ ß-\gamma+1\ \ ß-\gamma)$ | $6 \le \alpha \le 9$, $5 \le ß \le 8$, $2 \le \gamma \le 5$<br>$\alpha+\gamma \ge 11$, $9 \le ß+\gamma \le 10$, $2\alpha \ge 11+ß-\gamma$ |
| $P_8$ (124563) | $K_8(\alpha\ ß\ \gamma) = (\alpha-\gamma+1\ \ \alpha+ß-10\ \ ß-\gamma)$ | $6 \le \alpha \le 9$, $5 \le ß \le 8$, $1 \le \gamma \le 5$<br>$\alpha \ge ß+1$, $\alpha+ß \ge 10$, $10 \le \alpha+\gamma \le 11$, $9 \le ß+\gamma \le 11$ |
| $P_9$ (612453) | $K_9(\alpha\ ß\ \gamma) = (11-\alpha-\gamma\ \ \alpha+ß-9\ \ ß+\gamma-9)$ | $\alpha = 5$, $ß = 5$ $\gamma = 4$ |
| $P_1$ (1234) | $K_{22}(\alpha\ ß\ 0) = (\alpha-1\ \ ß\ \ \alpha-ß+1)$ | $6 \le \alpha \le 9$, $5 \le ß \le 9$, $\alpha \ge ß+1$, $2ß \ge \alpha+1$ |
| $P_2$ (1423) | $K_{24}(\alpha\ ß\ 0) = (ß\ \ \ 10-ß\ \ 2\alpha-10)$ | $5 \le \alpha \le 9$, $5 \le ß \le 9$, $\alpha+ß \le 11$ |
| $P_1$ (12) | $K_{31}(\alpha\ 0\ 0) = (\alpha-1\ \ 10-\alpha\ \ \ 0)$ | $6 \le \alpha \le 9$ |
| $P_2$ (21) | $K_{32}(\alpha\ 0\ 0) = (10-\alpha\ \ \alpha-1\ \ \ 0)$ | $0 \le \alpha \le 5$ |

Table 2. Some functions $K_i$ in $A_7$, w = 7, h = 3

| Permutation (*) | Associated $K_i$ (**α**) | Domain of existence |
|---|---|---|
| $P_1$ (123456) | $K_1$ (α ß γ) = (α-1  α+ß-9  ß+γ-9) | α≥ß+1, γ≥5 |
| $P_4$ (142563) | $K_4$ (α ß γ) = (10-γ  2α-10  ß-γ) | α≥ß+1, 9≤α+γ≤11, ß≥5, ß+γ≤9 |
| $P_5$ (145263) | $K_5$ (α ß γ) = (10-γ  2α-10  9-ß-γ) | α+ß≥10, ß≤4, 9≤α+γ≤11 |
| $P_6$ (412563) | $K_6$ (α ß γ) = (10-γ  α- γ-1  α+ß-9) | α≥ß+1, ß≥5, α+γ≤9 |
| $P_7$ (124536) | $K_7$ (α ß γ) = (α-1  α- γ+1  2ß-9) | α+γ≥11, 9≤ß+γ≤10 |
| $P_{10}$ (123465) | $K_{10}$ (α ß γ) = (ß  2α-10  ß+γ-9) | α≥γ+1, α≤ß+1, γ≥5 |
| $P_{11}$ (451623) | $K_{11}$ (α ß γ) = (10-γ  9-γ-ß  α-ß-1) | α≥5, α+ß≤9 |
| $P_{12}$ (145236) | $K_{12}$ (α ß γ) = (α-1  α- γ+1  9-ß-γ) | α+γ≥11, ß≤4 |
| $P_{13}$ (124365) | $K_{13}$ (α ß γ) = (ß  2α-10  ß-γ+1) | α≤ß+1, α+γ≥11, γ≤5 |
| $P_{22}$ (1234) | $K_{22}$ (α ß 0) = (α-1  ß  10-ß) | α≥ß+1, ß≥5 |
| $P_{23}$ (1243) | $K_{23}$ (α ß 0) = (ß  α-1  10-ß) | α≤ß+1, α+ß≤11 |
| $P_{25}$ (1324) | $K_{25}$ (α ß 0) = (α-1  10-ß  ß) | α+ß≤11, ß≤5 |
| $P_{26}$ (1342) | $K_{26}$ (α ß 0) = (10-ß  α-1  ß) | α≥ß+1, 9≤α+ß≤11 |
| $P_{27}$ (3142) | $K_{27}$ (α ß 0) = (10-ß  α-1  9- α) | α≥5, α+ß≤9 |
| $P_{31}$ (12) | $K_{31}$ (α 0 0) = (α-1  10-α  0) | 6≤α≤9 |
| $P_{32}$ (21) | $K_{32}$ (α 0 0) = (10-α  α-1  0) | 1≤α≤5 |

*Permutation of O (n'), excluding fixed nines.
Thus, $P_1$ (9  α  ß  γ-1  9-γ  9-ß  10-α), $P_{22}$ = (9  9  9  α  ß-1  9-ß  10-α),
$P_{31}$ = ( 9  9  9  9  9  α-1  10-α) and $P_6$ = ( 9  9-γ  α  ß  9-ß  10-α  γ-1)

## 4 Balance. Transformation constants

Let balance exist when there is an $\boldsymbol{\alpha}_E$ and a $K_i$ (**α**) such that $K_i$ ($\boldsymbol{\alpha}_E$) = $\boldsymbol{\alpha}_E$ or, equivalently, K ($n_E$) = $n_E$, p ($n_E$) = $\boldsymbol{\alpha}_E$

As a corollary $K^r(n_E) = n_E$, $K^r(\alpha_E) = \alpha_E$, $\forall r \geq 1$

Let $n_E$ be called **transformation constant** and **$\alpha_E$ parametric constant.** If all the numbers $n \in A_w$, by reiterating their transformation they become $n_E$ we will say that $n_E$ is a Kaprekar constant. This is equivalent to the existence of a constant and a single transformation tree in $A_w$.

A transformation constant is a single-link cycle, which is why some authors call the constants 1-cycle.

As an example of [4a] we have seen how 631764 transforms into itself, and the same is true for $\boldsymbol{\alpha} = p(n) = 632$, $K(632) = 632$. Therefore, $631764 = n_E$ and $632 = \boldsymbol{\alpha}_E$.

The constants in base 10 have been exhaustively studied by Prichett et al. (1981), who showed that in base 10 there are only two Kaprekar constants, 495 in $A_3$ and 6174 in $A_4$. They also determined 1-cycles for different w-digits.

Our objective in this section is not to make a detailed study of the constants, but to show how families of constants can be derived directly from the $K_i$ functions. In this regard we show two examples: the family $\boldsymbol{\alpha}_E = 6\ 3\ \overbrace{\ldots\ldots}^{h-2}\ 3\ 2$ in $w=2h$, $h \geq 2$ and the family $\boldsymbol{\alpha}_E = 8\ 6\ 4\ 3\ \overbrace{\ldots\ldots}^{h-4}\ 3\ 2$ in $w=2h+1$, $h \geq 7$. We also show two singular constants $\boldsymbol{\alpha}_E = 550$ in $A_6$ and $\boldsymbol{\alpha}_E = 5$ in $A_3$.

**<u>Constants $n_E = 63\ \overbrace{\ldots\ldots}^{h-2}\ 3176\ \overbrace{\ldots\ldots}^{h-2}\ 64$, $\boldsymbol{\alpha}_E = 63\ \overbrace{\ldots\ldots}^{h-2}\ 32$, $w = 2h$</u>**

Indeed $O_d(n_E) = 76\ \overbrace{\ldots\ldots}^{h-2}\ 6643\ \overbrace{\ldots\ldots}^{h-2}\ 31$ and $O_u(n_E) = 13\ \overbrace{\ldots\ldots}^{h-2}\ 3466\ \overbrace{\ldots\ldots}^{h-2}\ 67$

$n'_E = O_d(n_E) - O_u(n_E) = 63\ \overbrace{\ldots\ldots}^{h-2}\ 3176\ \overbrace{\ldots\ldots}^{h-2}\ 64 = n_E$

For this expression to be valid at least the two extreme digits and the two middle digits must be present. Such is the case with 6174, which is the renowned Kaprekar's constant. The associated parametric constants are $\boldsymbol{\alpha}_E = 63\ \overbrace{\ldots\ldots}^{h-2}\ 32$.

Such constants, which generalize Kaprekar's constant for w-digits, $w=\dot{2}$, derive from function [10] associated with permutation $P_3$ (h+1 1 h+2 … w 2 3 … h)

$K_3(\boldsymbol{\alpha}) = \boldsymbol{\alpha} \rightarrow \alpha^1 = 10-2\alpha^h$, $\alpha^2 = \alpha^1-\alpha^{h-1}$, $\alpha^3 = 9-\alpha^{h-1}-\alpha^{h-2} \ldots \alpha^h = \alpha^1-\alpha^2-1$

with the restrictions imposed in $K_3$. There is a unique solution which is $\boldsymbol{\alpha}_E = 6\ 3\ \overbrace{\ldots\ldots}^{h-2}\ 3\ 2$

For instance, for w = 8, h = 4

$K_3(\boldsymbol{\alpha}) = K_3(\alpha\ \beta\ \gamma\ \delta) = (10-2\delta\quad \alpha-\gamma\quad 9-\beta-\gamma\quad \alpha-\beta-1)$

$1 \leq 10-2\delta+\alpha-\gamma \leq 10$; $\beta+\gamma \leq 9$; $\alpha \geq \beta+1$, $\alpha+\gamma \geq 9$, $\alpha \leq 10+\gamma-2\delta$, $1 \leq \delta \leq 4$

Imposing the equality condition

$\alpha = 10\text{-}2\delta$, $ß = \alpha\text{-}\gamma$, $\gamma = 9\text{-}ß\text{-}\gamma$, $\delta = \alpha\text{-}ß\text{-}1$ there is only one solution $\alpha = 6$, $ß = 3$, $\gamma = 3$, $\delta = 2 \rightarrow \boldsymbol{\alpha}_E = 6332$ which satisfies the existence conditions and $f_1$ (6332) = 63317664 = $n_E$
In the case w = 4, h = 2, $P_3$ = (3142) and the associated function is
$K_3$ ($\alpha$ ß) = (10-2ß  2$\alpha$-10)  $\alpha \geq 5$, $\alpha + ß \leq 9$
$K_3$ ($\alpha$ ß) = ($\alpha$ ß) $\rightarrow \alpha = 10\text{-}2ß$, $ß = 2\alpha\text{-}10 \rightarrow \alpha = 6$, $ß = 2$
$f_1$ (6 2) = 6174, which is Kaprekar's constant.

This family of constants cannot exist if w = $\dot{2}$+1, since a 9 being present in [4b] completely changes function $K_3$. It is not valid for w = 2 either, since permutation $P_3$ does not make sense.

**Constants $n_E$ = 8643 $\overbrace{\ldots\ldots}^{h-4}$ 31976 $\overbrace{\ldots\ldots}^{h-4}$ 532, $\boldsymbol{\alpha}_E$ = 8643 $\overbrace{\ldots\ldots}^{h-4}$ 32, w=2h+1, h≥4**

They are derived from the function [10b] associated with the permutation $P_4$ (1  h+2  2  h+3 …h+s  3  4 … s  2h  2h+1  h), s = h+1, h≥4. De $K_4$ ($\boldsymbol{\alpha}$) = $\boldsymbol{\alpha}$ and taking into account the domain of existence of $P_4$, it results $\boldsymbol{\alpha}_E$.

**Constant $n_E$ = 549945 and $\boldsymbol{\alpha}_E$ = 550**

Derives from $K_{24}$ in $A_6$ (Table 1)
$\alpha = ß$, $ß = 10\text{-}ß$, $0 = 2\alpha\text{-}10 \rightarrow \alpha = 5$, $ß = 5$, $\gamma = 0 \leftrightarrow \boldsymbol{\alpha}_E = 550$
$f_2$ (550) = 549945 = $n_E$
This constant cannot be generalized.

**Constant $n_E$ = 495 and $\boldsymbol{\alpha}_E$ = 5**

This is the other Kaprekar's constant. If w = 3, $f_1$ ($\alpha$) = ($\alpha$-1  9  10-$\alpha$) and O (n') it allows for two permutations $P_{12}$ = (9 $\alpha$-1 10-$\alpha$) and $P_{21}$ = (9 10-$\alpha$ $\alpha$-1) as well as two associated functions $K_1$ ($\alpha$) = ($\alpha$-1) and $K_2$ ($\alpha$) = (10- $\alpha$). $\alpha_E$ = 5 results from $K_2$ ($\alpha$) = $\alpha$ and $f_1$ (5) = 495. It is thus another constant that cannot be generalized.

There are other functions of just one parameter. Such is the case of

- w = 2 with $f_1$ ($\alpha$) = ($\alpha$-1 10-$\alpha$) and two functions $K_i$,
  $K_1$ ($\alpha$) = (2$\alpha$-11) and $K_2$ ($\alpha$) = (11-2$\alpha$) which obviously cannot satisfy the condition $K_i$ ($\alpha$) = $\alpha$. In $A_2$ there cannot be any transformation constant.
- The same is true for functions
  $K_{31}$ ($\alpha$  0 $\overbrace{\ldots\ldots}^{h-1}$ 0) = ($\alpha$-1  10-$\alpha$  0 $\overbrace{\ldots\ldots}^{h-2}$ 0) and $K_{32}$ ($\alpha$  0 $\overbrace{\ldots\ldots}^{h-1}$ 0) =

  = (10-$\alpha$  $\alpha$-1  0 $\overbrace{\ldots\ldots}^{h-2}$ 0)

  where condition K ($\boldsymbol{\alpha}$) = $\boldsymbol{\alpha}$ leads to impossible results. In these cases there cannot be any transformation constant either.

## 5 Cycles: Terminology and functions $K^r$ ($\boldsymbol{\alpha}$)

There is a cycle of r links if and only of there are r parameters $\boldsymbol{\alpha}_{ci}$ and r functions $K_{ci}$ ($\boldsymbol{\alpha}$), such that $K_{ci}$ ($\boldsymbol{\alpha}_{ci}$) = $\boldsymbol{\alpha}_{ci+1}$ being $\boldsymbol{\alpha}_{ci} = \boldsymbol{\alpha}_{cj} \leftrightarrow i \equiv j \pmod r$.

This amounts to the existence of r operators $K^r_{ci}$ such that $K^r_{ci}(\boldsymbol{\alpha}_{ci}) = \boldsymbol{\alpha}_{ci}$, i = 1 …r
In terms of numbers, it demands the existence of r numbers
$n_{ci}$, $K(n_{ci}) = n_{ci+1}$, $n_{ci} = n_{cj} \leftrightarrow i \equiv j \pmod r$ or, equivalently, $K^r(n_{ci}) = n_{ci}$, i = 1…r.

Let us analyze the case $A_6$ for 6-digit numbers. The parametric transformation trees are shown in Graph 1.

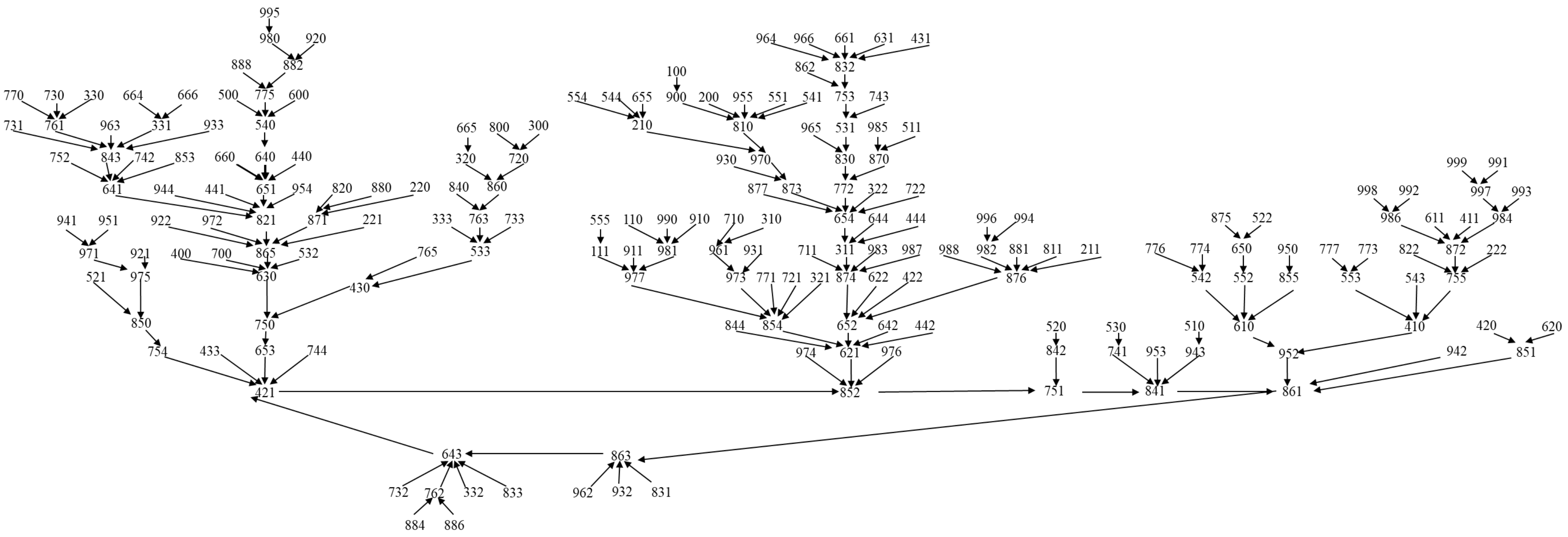

Tree A

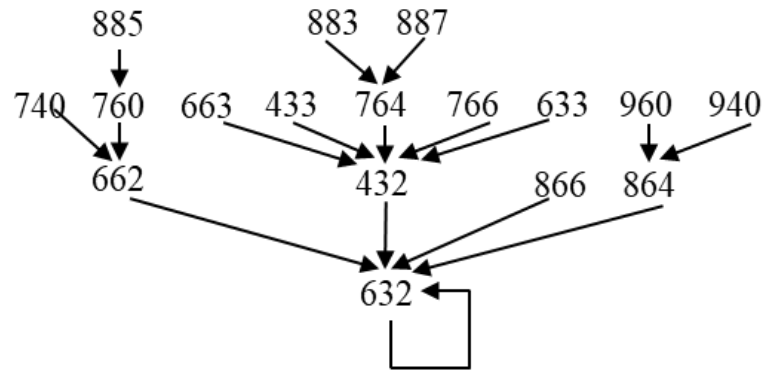

Tree B

Tree C

Graph 1. Parametric transformation trees in $A_6$. Each number is a class (α ß γ)

Main tree A, consisting of 201 parametric classes out of the 219 existing classes, is articulated on a 7-class cycle, which can be conventionally called
$861 = \boldsymbol{\alpha}_{c1}$, $863 = \boldsymbol{\alpha}_{c2}$, $643 = \boldsymbol{\alpha}_{c3}$, $421 = \boldsymbol{\alpha}_{c4}$, $852 = \boldsymbol{\alpha}_{c5}$, $751 = \boldsymbol{\alpha}_{c6}$, $841 = \boldsymbol{\alpha}_{c7}$

Each class becomes the following one with the functions in Table 1, indicated in the diagram

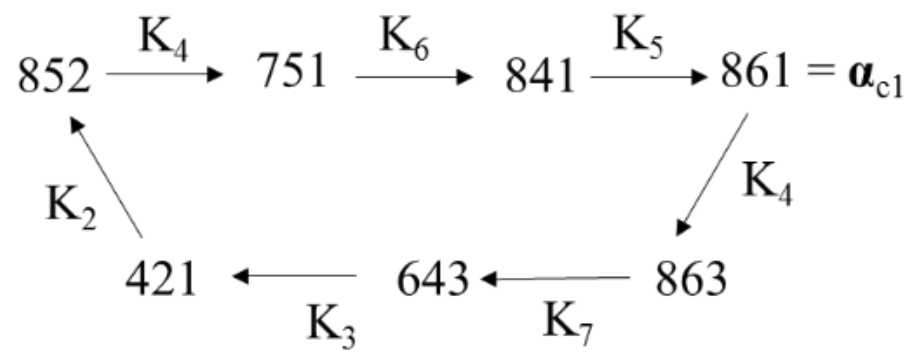

Naturally, there are 7 operators $K^r$, r = 1, 2…7, which transform parameters separated by r links. Thus,

$K^7_{c1}$ (861) = $K_5$ x $K_6$ x $K_4$ x $K_2$ x $K_3$ x $K_7$ x $K_4$ (861) =

= (20α+8ß-20γ-180   20α+8ß-20γ-178   24α-16 ß-24γ–71) = 861

Any operator $K^s_{ci}$ ($\boldsymbol{\alpha}_{ci}$) built with s > 7 will coincide with other $K^t_{ci}$ ($\boldsymbol{\alpha}_{ci}$) $1 \leq t \leq 7$ if and only if $s \equiv t \pmod 7$. **The existence of the cycle demands the coupling between functions $K_i$.** The image of one function must belong to the domain of existence of the following one.

In numerical terms the cycle consists of the following numbers

$n_{c1}$ = 840852, $n_{c2}$ = 860832, $n_{c3}$ = 862632, $n_{c4}$ = 642654, $n_{c5}$ = 420876, $n_{c6}$ = 851742, $n_{c7}$ = 750843, being p ($n_{ci}$) = $\boldsymbol{\alpha}_{ci}$.

In $A_2$ there is a 5-link cycle. In $A_3$ and $A_4$ there is only one tree articulated on a Kaprekar's constant. In $A_5$ there are two 4-link cycles and one 2-link cycle.

In $A_7$, consisting of 7-digit numbers, all 219 existing parametric classes (α ß γ) group with the following cycle formed by 8 links

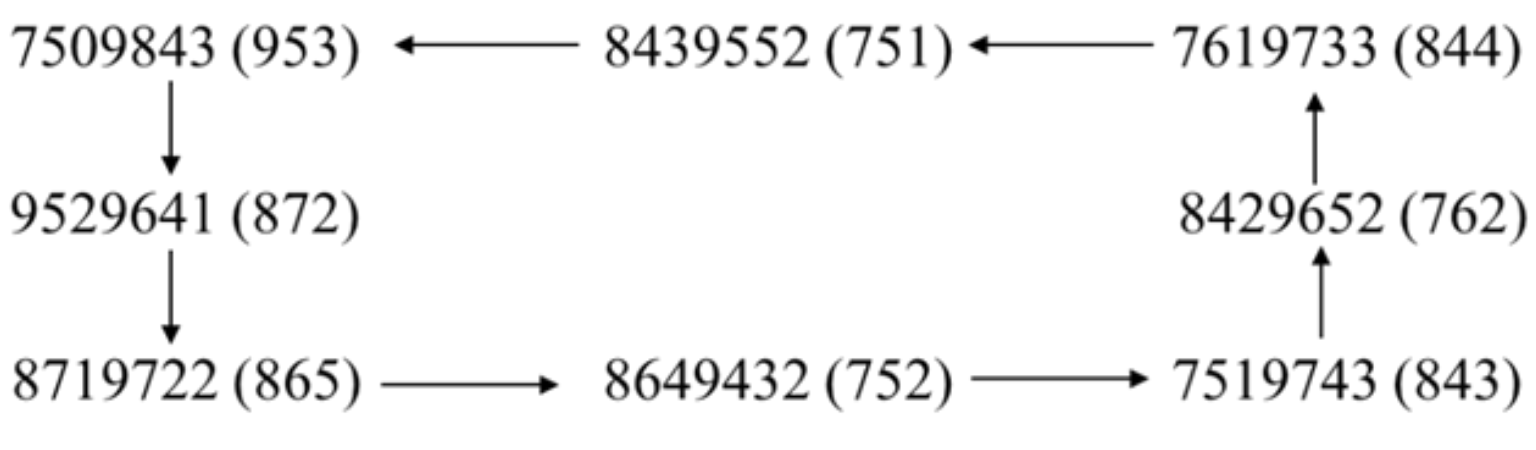

# 6 Symmetries

## 6.1 Equivalence relations

Let two numbers m and n be equivalent $R_r$ of order r if and only if
$m \, R_r \, n \leftrightarrow K^r(m) = K^r(n) \qquad r \geq 0$ [20]
Since $K^r(m) = K^{r-1}[K(m)]$ the following recurrence equations result

$m\ R_r\ n \leftrightarrow K(m)\ R_{r-1}\ K(n) \quad r \geq 1$ [21]

$m\ R_r\ n \leftrightarrow K^s(m)\ R_{r-s}\ K^s(n) \quad s \leq r$ [21a]

Note that if two numbers are equivalent, they are so for any equivalence of a higher order

$m\ R_r\ n \rightarrow m\ R_{r+s}\ n \quad s \geq 0$ [22]

It is necessary to distinguish **new equivalences $R_r$**, those which are not of a lower order

$m\ \mathbf{R_r}\ n \rightarrow m\ \not R_{r-1}\ n \quad r > 0$ [23]

$m\ \mathbf{R_r}\ n \rightarrow m\ \not R_{r-s}\ n \quad s = 1,2\ldots r$ [23a]

from **old** equivalences ${}_aR_r$, those that are also of a lower order

$m\ ({}_aR_r)\ n \rightarrow m\ R_{r-1}\ n \quad r > 0$ [24]

The notation $R_r$ does not specify whether an equivalence is new or old

Equivalence of an order rare binary equivalence relations that have the following features

reflexivity $m\ R_r\ m \quad \forall\ m \in A_W$

symmetry $m\ R_r\ n \leftrightarrow n\ R_r\ m \quad \forall\ m,n$

transitivity $m\ R_r\ n$ and $n\ R_r\ q \rightarrow m\ R_r\ q \quad \forall\ m,n,q$

They allow for numbers to be grouped in equivalence classes of order r

$C_r = \{m,n \in A_W;\ m\ R_r\ n\}$ [25]

These equivalence classes are pairwise disjoint and their intersection is set $A_W$. generating a partition in $A_W$.

Because two numbers that are equivalent are so for any equivalence of a higher order [22], classes $C_r$ will include both the new and old equivalences. This hierarchical and inclusive feature of binary relations $R_r$ is essential to better understand the convergence of this process.

$C_r = \{m,n \in A_W,\ m\ \mathbf{R_r}\ n\},\ {}_aC_r = \{\ p,q \in A_W,\ p\ R_{r-1}\ q\} \quad C_r = \mathbf{C_r} \cup {}_aC_r$ [26]

The transitive feature and relation [22] yield

$m\ R_r\ n$ and $n\ R_{r+s}\ q \rightarrow m\ R_{r+s}\ q \quad s \geq 0$ [27]

In our methodological approach to the analysis of the Kaprekar process, it is essential to extend the concept of equivalence $R_r$ defined between numbers to sets of numbers with images common, to equivalence classes.

Considering the last relationship, if

$m \in C_r^i,\ p \in C_s^j$ and $m\ R_t\ p,\ t = \max(r,s) + u,\ u \geq 0$ then if

$n \in C_r^i$ and $q \in C_s^j$ results $n\ R_t\ q$. This is due to the univocal character of the transformation. All numbers belonging to a class $C_v$ have the same image in the v-th transformation. From here on, any subsequent transformation will have the same image. Consequently, the element that best represents the transformation is not the number, but the equivalence class. Thus, we can agree to define the equivalence relation between classes according to

$m \in C_r^i,\ n \in C_s^j$ and $m\ R_t\ p,\ t = \max(r,s) + u,\ u \geq 0 \leftrightarrow C_r^i\ R_t\ C_s^j$ [10]

The direction of the arrow to the left must be understood in the sense that if two classes $C_r^i$ y $C_s^j$ have a relation $R_t$ then all numbers m and n belonging to those classes have a relation $R_t$. In particular it turns out

$m \in C_r^i$, $n \in C_r^j$, $m\ R_{r+s}\ n \leftrightarrow C_r^i\ R_{r+s}\ C_r^j\ \ s \geq 0$ [28a]

From relation [7] it follows that the $C_1$ equivalence classes are defined by the parametric classes (we will return to this topic in paragraph 6.3). $\boldsymbol{\alpha}$ represents the parameters of a number $\boldsymbol{\alpha} = (\alpha^1\ \alpha^2 \dots \alpha^h) = p(n)$ and also the set of numbers that has the same parameters. As what characterizes a class $C_1$ are precisely these parameters, we can agree to specify each $C_1$ by its corresponding parameters, being able to write $C_1^i = \boldsymbol{\alpha}_i$, if in the context there is no room for confusion. With this convention $\boldsymbol{\alpha}$ has a double meaning: parameters of a number and set of numbers with the parameters that are specified.

By extension, we will say that two parametric classes $\boldsymbol{\alpha}_1$ and $\boldsymbol{\alpha}_2$ are equivalent $R_r$ if and only if the numbers with the respective parameters are equivalent

$\boldsymbol{\alpha}_1\ R_r\ \boldsymbol{\alpha}_2 \leftrightarrow m\ R_r\ n,\ \ p(m) = \boldsymbol{\alpha}_1,\ \ p(n) = \boldsymbol{\alpha}_2$ [28b]

The expression [28b] allows us to translate some of the previous relationships into a parametric class language. For example from [28a] it results

$\alpha_1 \subset C_r^i,\ \ \alpha_2\ C_r^j,\ \ \alpha_1\ R_{r+s}\ \alpha_2\ \rightarrow C_r^j\ R_{r+s} C_r^j,\ \ s \geq 0$ [29]

Going back to the expression [7] and to the definition [20] it results

$m\ R_r\ n \leftrightarrow p\,[K^{r-1}(m)] = p\,[K^{r-1}(n)] \qquad r \geq 1$ [30]

To any equivalence $R_r$ can be associated a function or operator $e_r$ such that

$\boldsymbol{\alpha}_1\ R_r\ \boldsymbol{\alpha}_2 \leftrightarrow e_r(\boldsymbol{\alpha}_1) = \boldsymbol{\alpha}_2$ [31]

with the domain of existence of $R_r$.

It is a kind of representation that makes calculation easier.

### 6.2 <u>Product of equivalences</u>

Applying an equivalence onto another one normally yields another equivalence. By using the usual representation of the product of functions

$e^i\,[e^j(\boldsymbol{\alpha})] = (e^i \times e^j)(\boldsymbol{\alpha}) = e^k(\boldsymbol{\alpha})$

the domains of existence will have to be compatible. Thus, taking equivalences that will be later deduced

$e_2^2(\alpha\ ß) = (10\text{-}\alpha\ \ ß),\ \alpha+ß \leq 10;\ e_2^3(\alpha\ ß) = (\alpha\ \ 10\text{-}ß),\ \alpha+ß \geq 10$

$e_2^4(\alpha\ ß) = (e_2^3 \times e_2^2)(\alpha\ ß) = (10\text{-}\alpha\ 10\text{-}ß),\ \alpha = ß$

However,

$(e_2^2 \times e_2^4)(\alpha\ ß) = (\alpha\ 10\text{-}ß) = e_2^3$, where domains $\alpha = ß$ and $\alpha+ß \leq 10$ lead to $\alpha+ß \geq 10$.

Moreover, the product is conditioned by restriction [3]. Therefore, if $\boldsymbol{\alpha} = (\alpha\ ß\ \gamma)$, $\alpha \geq ß \geq \gamma$

$e_2^2(\boldsymbol{\alpha}) = (\alpha\ 9\text{-}ß\ \gamma),\ \ \alpha+ß \geq 9,\ ß+\gamma \leq 9;\ e_2^3(\boldsymbol{\alpha}) = (\alpha\ \ ß\ \ 10\text{-}\gamma),\ ß+\gamma \geq 10$

$e_2^6(\boldsymbol{\alpha}) = (e_2^3 \times e_2^2)(\boldsymbol{\alpha}) = (\alpha\ \ 9\text{-}ß\ \ 10\text{-}\gamma)$ not valid because it demands $\gamma \geq ß+1$, which contradicts [3].

And the other way around, the product of non-allowed equivalences can generate an allowed one. Thus,

$e_2^7(\boldsymbol{\alpha}) = (10\text{-}\alpha\ \ 9\text{-}ß\ \ 10\text{-}\gamma)$, not valid

$(e_2^6 \times e_2^7)(\boldsymbol{\alpha}) = (10\text{-}\alpha\ \ ß\ \ \gamma) = e_2^1(\boldsymbol{\alpha})$ valid if $\alpha+ß \leq 10$

In sum, the product of equivalences demands that the resulting equivalence's domain of existence be determined.

### 6.3 Symmetries $R_0$ and $R_1$

Symmetry $R_0$ has been defined conventionally as the no-transformation: $m\ R_0\ n \leftrightarrow K^0(m) = K^0(n) \leftrightarrow m = n$. Each equivalence class $C_0$ consists of a single number. There are as many classes $C_0$ as there are numbers belonging to $A_W$.

Due to [7], symmetry $R_1$ is the one containing those numbers with the same parameters, which return the same number after one transformation. Therefore, parametric classes, understood as the set of numbers with the same parameters, match with equivalence classes $C_1$. Graph 1 shows all these classes and their transformations in the case of $A_6$.
As these classes uniquely determine the image of any number belonging to them, they are an essential element when understanding Kaprekar's routine.

### 6.4 Symmetries $R_2$

As a result of $R_r$ definition [20]
$m\ R_2\ n \leftrightarrow K^2(m) = K^2(n)$ and because of [30]
$m\ R_2\ n \leftrightarrow p\,[K(m)] = p\,[K(n)]$ [32]

#### 6.4.1 Permutations of the same sequence

If two numbers m and n, $p(m) = \boldsymbol{\alpha}_1$, $p(n) = \boldsymbol{\alpha}_2$ return images with transpositions of digits compatible with belonging to $B_w$, the ordered sequences will be identical $O_d(m') = O_d(n')$ and their parameters will match $\boldsymbol{\alpha}'_1 = \boldsymbol{\alpha}'_2$, as will their numeric images $m'' = f(\boldsymbol{\alpha}'_1) = f(\boldsymbol{\alpha}'_2) = n''$, thus verifying the condition of $m\ R_2\ n$.
Not just any transposition is valid. Only those that respect the conditions [4a.1] or [4b.1] of belonging to $B_w$ are valid.

**Equivalences $R_2$ based on [4a] and [4b]**

$m' = f_1(\boldsymbol{\alpha}_1) = (\alpha_1^1\ \ \alpha_1^2 \ldots \alpha_1^{h-1}\ \ \alpha_1^h\text{-}1\ \ \ 9\text{-}\alpha_1^h \ldots 9\text{-}\alpha_1^2\ \ \ 10\text{-}\alpha_1^1)$, $p(m) = \boldsymbol{\alpha}_1 = (\alpha_1^1\ \alpha_1^2 \ldots \alpha_1^h)$

If $p(n) = \boldsymbol{\alpha}_2 = (10\text{-}\alpha_1^1\ \alpha_1^2 \ldots \alpha_1^h)$ $\quad n' = f_1(\boldsymbol{\alpha}_2) = (10\text{-}\alpha_1^1 \ldots \alpha_1^2 \ldots 9\text{-}\alpha_1^2\ \alpha_1^1)$

$n' \in B_w$ because it meets the conditions [4a.1] or [4b.1] and the sequence of n' has the same digits as that of m', the ends being translocated. Consequently,

$O_d(m') = O_d(n')$ and $p(m') = p(n') \leftrightarrow \boldsymbol{\alpha}_1\ R_2\ \boldsymbol{\alpha}_2$

We used $f_1$ corresponding to $w = \dot{2}$, [4a]. The result for $w = \dot{2}+1$ does not change, since including a central 9 in the sequences of m' and n' does not modify the reasoning. Likewise, the symmetric digits can be transposed in relation to the middle ones, $\alpha^s$ and $9\text{-}\alpha^s$, $s = 2, 3 \ldots h\text{-}1$, or the middle digits themselves, $\alpha^h\text{-}1$ and $9\text{-}\alpha^h$, returning the following binary equivalences by simple transposition, where $\boldsymbol{\alpha}_1 = (\alpha^1\ \alpha^2 \ldots \alpha^h)$.

- **Equivalences $\alpha_1\ R_2^{11}\alpha_2$, $\alpha_2 = (10\text{-}\alpha^1\ \alpha^2 \ldots \alpha^h)$, $\alpha^1+\alpha^2 \leq 10$**

- **Equivalences $\alpha_1\ R_2^{12}\alpha_2$, $\alpha_2 = (\alpha^1\ \alpha^2 \ldots \alpha^{s-1}\ 9\text{-}\alpha^s\ \ \alpha^{s+1} \ldots \alpha^h)$, s = 2, …h-1**
  **$\alpha^{s-1} + \alpha^s \geq 9$, $\alpha^s + \alpha^{s+1} \leq 9$**
  For example: w = 11, $\boldsymbol{\alpha}_1$ = 76632, $\boldsymbol{\alpha}_2$ = 76332
  $f_1$ ($\boldsymbol{\alpha}_1$) = 76631976333, $f_1$ ($\boldsymbol{\alpha}_2$) = 76331976633 being
  $f_1$ ($\boldsymbol{\alpha}_2$) the permutation P = (1 2 9 4 5 6 7 8 3 10 11) of $f_1$ ($\boldsymbol{\alpha}_1$) and K ($\boldsymbol{\alpha}_1$) = K ($\boldsymbol{\alpha}_2$) = 84433 consistent with [32]. Besides, $f_1$ (8 4 4 3 3) = ( 8 4 4 3 2 9 6 6 5 5 2) =
  = K (7 6 6 3 1 9 7 6 3 3) = K ( 7 6 3 3 1 9 7 6 6 3 3).

- **Equivalences $\alpha_1\ R_2^{13}\alpha_2$, $\alpha_2 = (\alpha^1\ \alpha^2 \ldots \alpha^{h-1}\ \ 10\text{-}\alpha^h)$, $\alpha^{h-1} + \alpha^h \geq 10$**
  However, multiple transpositions are limited by [3]:

- **Non-equivalence $\alpha_1\ n_2^1\ \alpha_2$**
  $\boldsymbol{\alpha}_2 = (\alpha^1 \ldots \alpha^{s-1}\ 9\text{-}\alpha^s\ \ \alpha^{s+1} \ldots \alpha^{s+r-1}\ 9\text{-}\alpha^{s+r}\ \alpha^{s+r+1} \ldots \alpha^h)$ impossible, $\boldsymbol{\alpha}_1\ \not R_2\ \boldsymbol{\alpha}_2$
  $\alpha^{s+r-1} \leq 9\text{-}\ \alpha^s$ and $\alpha^{s+r-1} \geq 9\text{-}\alpha^{s+r}$
  $\alpha^{s+r} \geq 9\text{-}\ \alpha^{s+r-1} \geq \alpha^s \rightarrow \alpha^{s+r} = \alpha^s = \alpha^{s+i}$ , i = 0,1,…r
  $9\text{-}\ \alpha^s \geq \alpha^s \rightarrow \alpha^s \leq 4 \rightarrow \alpha^{s+r-1} = \alpha^{s+r} \leq 4$
  $\alpha^{s+r-1} \geq 9\text{-}4 \geq 5$, which contradicts the previous expression. Therefore, $\boldsymbol{\alpha}_1\ \not R_2\ \boldsymbol{\alpha}_2 \leftrightarrow \boldsymbol{\alpha}_1\ n_2^1\ \boldsymbol{\alpha}_2$

- **Non-equivalence $\alpha_1\ n_2^2\ \alpha_2$**
  $\boldsymbol{\alpha}_2 = (\alpha^1 \ldots \alpha^{s-1}\ 9\text{-}\alpha^s\ \ \alpha^{s+1} \ldots 10\text{-}\alpha^h)$ impossible, $\boldsymbol{\alpha}_1\ \not R_2\ \boldsymbol{\alpha}_2$
  $\alpha^{h-1} \leq 9\text{-}\ \alpha^s$ and $\alpha^{h-1} \geq 10\text{-}\alpha^h$
  $\alpha^h \geq 10\text{-}\ \alpha^{h-1} \geq 10\text{-}(9\text{-}\ \alpha^s) = 1+\ \alpha^s$, impossible, since $\alpha^h \leq \alpha^s$
  $\boldsymbol{\alpha}_1\ \not R_2\ \boldsymbol{\alpha}_2 \leftrightarrow\ \boldsymbol{\alpha}_1\ n_2^2\ \boldsymbol{\alpha}_2$

- **Equivalence $\alpha_1\ R_2^{14}\ \alpha_2$,** $\boldsymbol{\alpha}_2 = (10\text{-}\alpha^1\ \alpha^2 \ldots \alpha^{h-1}\ \ 10\text{-}\alpha^h) \rightarrow \boldsymbol{\alpha}_2 = (\overbrace{5 \ldots 5}^{h}) = \boldsymbol{\alpha}_1$

  It is a specific identity solution
  $\alpha^{h-1} \leq 10\text{-}\ \alpha^l$ y $\alpha^{h-1} \geq 10\text{-}\ \alpha^h \rightarrow \alpha^h \geq \alpha^1 \rightarrow \alpha^h = \alpha^1 = \alpha^s$, s = 1,...h

  $$\left.\begin{array}{l} 10\text{-}\ \alpha^1 \geq \alpha^1 \rightarrow \alpha^1 \leq 5 \rightarrow \alpha^{h-1} = \alpha^h \leq 5 \\ \alpha^{h-1} \geq 10\text{-}\ \alpha^h \geq 10\text{-}5 = 5 \end{array}\right\}\ \alpha^{h-1} = 5\ \rightarrow\ \alpha^s = 5,\ s = 1 \ldots h$$

- **Equivalence of $\alpha\ R_2^{15}\ \alpha_2$, $\alpha_2 = (\alpha^1 \ldots \alpha^{s-1}\ 9\text{-}\alpha^s \ldots 9\text{-}\alpha^{s+r}\ \ \alpha^{s+r+1}\ \ \alpha^h)$**

  According to [3], $9\text{-}\ \alpha^s \geq 9\text{-}\alpha^{s+r} \rightarrow \alpha^{s+r} = \alpha^s$
  Domain of existence:
  $\alpha^s \geq 9\text{-}\ \alpha^{s-1}$, $\alpha^{s+r+1} \leq 9\text{-}\ \alpha^{s+r}$, $\alpha^{s+i} = \alpha^s$, i = 0, 1 …r+1, $2 \leq s < h$
  For example: w = 10
  $\boldsymbol{\alpha}_1$ = 75552 and $\boldsymbol{\alpha}_2$ = 74442, $f_1$ ($\boldsymbol{\alpha}_1$) = 7555174443, $f_1$ ($\boldsymbol{\alpha}_2$) = 7444175553, where $f_1$ ($\boldsymbol{\alpha}_2$) is the permutation P = (1 9 8 7 5 6 4 3 2 10) of $f_1$ ($\boldsymbol{\alpha}_1$) and K ($\boldsymbol{\alpha}_1$) = K ($\boldsymbol{\alpha}_2$) = = 64111 according to [32]

- **Equivalence $\alpha_1\ R_2^{16}\ \alpha_2$, $\alpha_2 = (10\text{-}\alpha^1\ \ \alpha^2\ \ 9\text{-}\alpha^3 \ldots 9\text{-}\alpha^r\ \ \alpha^{r+1}\ \ \alpha^h)$**
  $9\text{-}\alpha^3 \leq \alpha^2 \leq 10\text{-}\ \alpha^1 \rightarrow \alpha^2 \leq \alpha^3 \geq \alpha^1\text{-}1 \rightarrow \alpha^1 = \alpha^2 = \alpha^3 = 5$

Domain of existence:
$\alpha^s = 5$, s = 1, 2…r, $\alpha^{r+1} \leq 9- \alpha^r$
For example, $\boldsymbol{\alpha}_1$ = 555532, $\boldsymbol{\alpha}_2$ = 554432, w = 12
$f_l$ ($\boldsymbol{\alpha}_1$) = 555531764445, $f_l$ ($\boldsymbol{\alpha}_2$) = 554431765545
where $f_l$ ($\boldsymbol{\alpha}_2$) is a permutation P = (1 2 10 9 5 6 7 8 4 3 11 12) of $f_l$ ($\boldsymbol{\alpha}_1$) and K ($\boldsymbol{\alpha}_1$) = K ($\boldsymbol{\alpha}_2$) = 631110

- **Equivalence $\boldsymbol{\alpha}_1$ $R_2^{17}$ $\boldsymbol{\alpha}_2$, $\boldsymbol{\alpha}_2$ = (10-$\alpha^1$ 9-$\alpha^2$….9-$\alpha^r$ $\alpha^{r+1}$ $\alpha^h$)**
  $\alpha^1$-1 $\leq \alpha^2 \leq \alpha^1$, $\alpha^s = \alpha^2$, s = 2, 3,…r
  For example: w = 9, $\boldsymbol{\alpha}_1$ = 8772, $\boldsymbol{\alpha}_2$ = 2222
  $f_l$ ($\boldsymbol{\alpha}_1$) = 877197222, $f_l$ ($\boldsymbol{\alpha}_2$) = 222197778 where $f_l$ ($\boldsymbol{\alpha}_2$) is the permutation P = (9 8 74 5 6 3 2 1) of $f_l$ ($\boldsymbol{\alpha}_1$) and K ($\boldsymbol{\alpha}_1$) = K ($\boldsymbol{\alpha}_2$) = 8655

**<u>$R_2$ based on [5a] or [5b]</u>**

$\boldsymbol{\alpha} = (\alpha^1 \ldots \alpha^r\ \underbrace{0 \ldots 0}_{h-r})$ and functions $f_2$ ($\boldsymbol{\alpha}$) [5a] or [5b] operate

If we consider $\alpha^r$ as $\alpha^h$ the previous equivalences are valid. For example,

$\boldsymbol{\alpha}\ R_2^2\ (\alpha^1 \ldots \alpha^{r-1}\ 10\text{-}\alpha^r\ \underbrace{0 \ldots 0}_{h-r})$ $\boldsymbol{\alpha}_1$ = 97400 $R_2$ 97600 = $\boldsymbol{\alpha}_2$, since K ($\boldsymbol{\alpha}_1$) = K ($\boldsymbol{\alpha}_2$) =

= 87642

**<u>$R_2$ based on [6]</u>**

$\boldsymbol{\alpha} = (\alpha^1\ \underbrace{0 \ldots 0}_{h-1})$ and the function $f_3$ ($\boldsymbol{\alpha}$) = $(\alpha^1\text{-}1\ \underbrace{9 \ldots 9}_{w-2}\ 10\text{-}\alpha^1)$ operates

Only the transposition of extreme digits is possible. That is, if

$\boldsymbol{\alpha}_1 = (\alpha^1\ \underbrace{0 \ldots 0}_{h-1})$ and $\boldsymbol{\alpha}_2 = (11\text{-}\alpha^1\ \underbrace{0 \ldots 0}_{h-1})$ then

$f_3(\boldsymbol{\alpha}_1) = (\alpha^1\text{-}1\ \underbrace{9 \ldots 9}_{w-2}\ 10\text{-}\alpha^1)$, $f_3(\boldsymbol{\alpha}_2) = (10\text{-}\alpha^1\ \underbrace{9 \ldots 9}_{w-2}\ \alpha^1\text{-}1)$

which is a permutation of the same digits compatible with [6.1], so
$O_d$ [$f_3$ ($\boldsymbol{\alpha}_1$)] = $O_d$ [$f_3$ ($\boldsymbol{\alpha}_2$)] ↔ $\boldsymbol{\alpha}'_1$ = K ($\boldsymbol{\alpha}_1$) = $\boldsymbol{\alpha}'_2$ = K ($\boldsymbol{\alpha}_2$) ↔ $\boldsymbol{\alpha}_1$ $R_2$ $\boldsymbol{\alpha}_2$

- **Equivalence $(\alpha\ \underbrace{0 \ldots 0}_{h-1})\ R_2^3\ (11\text{-}\alpha\ \underbrace{0 \ldots 0}_{h-1})$, $\alpha \geq 2$**

  Example: $\alpha_1$ = 300, $\alpha_2$ = 800, $\alpha_1$ $R_2$ $\alpha_2$ since K ($\alpha_1$) = K ($\alpha_2$) = 720

### 6.4.2 <u>Permutation of different sequences</u>

There are equivalences where the digits from the sequence do not necessarily remain.

- **Equivalence $\boldsymbol{\alpha}_1$ $R_2^4$ $\boldsymbol{\alpha}_2$,**
  **$\boldsymbol{\alpha}_1$=($\alpha^1$ $\alpha^2$ … $\alpha^h$), $\boldsymbol{\alpha}_2$=(10-$\alpha^h$ 9-$\alpha^{h-1}$ 9-$\alpha^{h-2}$ … 9-$\alpha^3$ 9-$\alpha^2$ 10-$\alpha^1$),**
  $\alpha^1 \geq \alpha^2$+1, w = $\dot{2}$, h = w/2
  Under [8a] $K_1$ ($\boldsymbol{\alpha}_1$) = (2 $\alpha_1^1$-10 2 $\alpha_1^2$-9 … 2 $\alpha_1^{h-1}$-9 2 $\alpha_1^h$-10)

Under [9a] $K_2$ ($\boldsymbol{\alpha}_2$) = (10-2 $\alpha_2^h$ 9-2 $\alpha_2^{h-1}$ … 9-2 $\alpha_2^2$ 10-2 $\alpha_2^1$)

If $K_1$ ($\boldsymbol{\alpha}_1$) = $K_2$ ($\boldsymbol{\alpha}_2$) → $\alpha_1^1 + \alpha_2^h = 10$, $\alpha_1^{s+1} + \alpha_2^{h-s} = 9$, s = 1, 2, … h-2

Note that $\boldsymbol{\alpha}_1$ + $O_u$ ($\boldsymbol{\alpha}_2$) = (10 9 $\underbrace{\ldots\ldots}_{h-2}$ 9 10)

For example;

m = 8178382562, p (m) = $\boldsymbol{\alpha}_1$ = 76641, m' = $f_1$ ($\boldsymbol{\alpha}_1$) = 7664085333

n = 4774473809, p (n) = $\boldsymbol{\alpha}_2$ = 95333, n' = $f_1$ ($\boldsymbol{\alpha}_2$) = 9533266641

$\boldsymbol{\alpha}_1$ and $\boldsymbol{\alpha}_2$ verify the equivalence, and $f_1$ ($\boldsymbol{\alpha}$) corresponds to [4a]. m' and n' do not keep the same digits, but p (m') = p (n') = 84331 and $K^2$ (m) = $K^2$ (n) = 8433086652. Thus answering the second question in the introduction.

This equivalence, contrary to what happens with the previous ones, is not valid for odd and w. In the previous example w=10. Let w = 11 and the same parameters p (m) = $\boldsymbol{\alpha}_1$ = 76641, p (n) = $\boldsymbol{\alpha}_2$ = 95333

Now the function $f_1$ that must be used is the one corresponding to [4b]

m' = $f_1$ ($\boldsymbol{\alpha}_1$) = 76640985333, n'=$f_1$ ($\boldsymbol{\alpha}_2$) = 95332966641

which are similar to the previous ones but with a 9 in the middle position

p (m') = 95432, p (n') = 87332

Adding a 9 in the sequence of digits has led the symmetry to break.

This equivalence may generate permutations of digits from the same sequence. Such is the case of $\boldsymbol{\alpha}_1$ = 87662 and $\boldsymbol{\alpha}_2$ = 83322, $\boldsymbol{\alpha}_1$ $R_2^4$ $\boldsymbol{\alpha}_2$

m' = $f_1$ ($\boldsymbol{\alpha}_1$) = 8766173322, n'=$f_1$ ($\boldsymbol{\alpha}_2$) = 8332177662

being $O_d$ (m') = $O_d$ (n'). It is just a numeric coincidence, even if for these parameters $R_2^4$ does not match any of the previous equivalences. **$R_2^4$ is a very common equivalence** given the few restrictions of its domain of existence.

The product of equivalences $R_2^1$ x $R_2^4$ creates a family of derivative equivalences. Table 3 shows the equivalences $R_2$ common in $A_6$ and $A_7$, and Table 4 shows derivative equivalences in $A_6$. The superscripts in the equivalences have been simplified by using correlative natural numbers. $n_2^6$ corresponds to a double translocation, and $n_2^7$ to a triple one, but none of them is an equivalence.

Table 3. Equivalences $R_2$ with three parameters, $\boldsymbol{\alpha}$ = (α ß γ).

| Equivalences $R_2$ common to $A_6$ and $A_7$ | | |
|---|---|---|
| Group | $e_2^i$ ($\boldsymbol{\alpha}_1$) = $\boldsymbol{\alpha}_2$ | Domain of existence |
| Set I Based on [4a] or [4b] $\boldsymbol{\alpha}$ = (α ß γ) 1≤α≤9 1≤ß≤9 1≤γ≤9 | $e_2^0$ ($\boldsymbol{\alpha}$) = (α ß γ) | 1≤α≤9, 1≤ß≤9, 1≤γ≤9 |
| | $e_2^1$ ($\boldsymbol{\alpha}$) = (10-α ß γ) | α+ß ≤ 10 |
| | $e_2^2$ ($\boldsymbol{\alpha}$) = (α 9-ß γ) | α+ß ≥ 9, ß+γ ≤ 9 |
| | $e_2^3$ ($\boldsymbol{\alpha}$) = (α ß 10-γ) | ß+γ ≥10 |
| | $e_2^4$ ($\boldsymbol{\alpha}$) = (10-α 9-ß γ) | α ≤ ß+1, ß+γ ≤ 9 |
| | $e_2^5$ ($\boldsymbol{\alpha}$) = (10-α ß 10-γ) | α = ß = γ = 5 |
| | $n_2^6$ ($\boldsymbol{\alpha}$) = (α 9-ß 10-γ) | α+ß ≥ 9, γ ≥ ß+1→ Impossible |
| | $n_2^7$ ($\boldsymbol{\alpha}$) = (10-α 9-ß 10-γ) | α≤ß+1, γ ≥ ß+1→ Impossible |
| Group II Based on [5a] or [5b] $\boldsymbol{\alpha}$ = (α ß 0) | $e_2^0$ ($\boldsymbol{\alpha}$) = (α ß 0) | 1≤α≤9, 1≤ß≤9, γ=0 |
| | $e_2^8$ ($\boldsymbol{\alpha}$) = (10-α ß 0) | α+ß ≤ 10 |
| | $e_2^9$ ($\boldsymbol{\alpha}$) = (α 10-ß 0) | α+ß ≥ 10 |
| | $e_2^{10}$ ($\boldsymbol{\alpha}$) = (10-α 10-ß 0) | α = ß |
| Group III Based on [6] $\boldsymbol{\alpha}$ = (α 0 0) | $e_2^0$ ($\boldsymbol{\alpha}$) = (α 0 0) | 1≤α≤9 |
| | $e_2^{11}$ ($\boldsymbol{\alpha}$) = (11-α 0 0) | 2≤α≤9 |
| | Equivalences $R_2$ in $A_6$ | |
| Group IV $\boldsymbol{\alpha}$ = (α ß γ) | $e_2^0$ ($\boldsymbol{\alpha}$) = (α ß γ) | 1≤α≤9, 1≤ß≤9, 1≤γ≤9 |
| | $e_2^{12}$ ($\boldsymbol{\alpha}$) = (10- γ 9-ß 10-α) | 1≤α≤9, 1≤ß≤8, 1≤γ≤8, α ≥ ß+1 |

Table 4. Derivative equivalences in $A_6$, $\boldsymbol{\alpha}$ = (α ß γ)

| Equivalence | Domain of existence | Example |
|---|---|---|
| $(e_2^1 \text{ x } e_2^{12})(\boldsymbol{\alpha})$ = (γ 9-ß 10-α) | $\alpha \geq ß+1, ß+\gamma \geq 9$ | 963 $R_2$ 331 |
| $(e_2^2 \text{ x } e_2^{12})(\boldsymbol{\alpha})$ = (10-γ ß 10-α) | $\alpha+ß \geq 10, ß+\gamma \leq 10$ | 963 $R_2$ 761 |
| $(e_2^3 \text{ x } e_2^{12})(\boldsymbol{\alpha})$ = (10-γ 9-ß α) | $\alpha+ß \leq 9$ | 221 $R_2$ 972 |
| $(e_2^4 \text{ x } e_2^{12})(\boldsymbol{\alpha})$ = (γ ß 10-α) | $\alpha+ß \geq 10, ß = \gamma$ | 933 $R_2$ 331 |
| $(n_2^6 \text{ x } e_2^{12})(\boldsymbol{\alpha})$ = (10-γ ß α) | $\alpha = ß, ß+\gamma \leq 10$ | 772 $R_2$ 877 |
| $(n_2^7 \text{ x } e_2^{12})(\boldsymbol{\alpha})$ = (γ ß α) | $\alpha = ß = \gamma$ | 333 $R_2$ 333 |

Some of these derivative equivalences are not reciprocal, but

$\boldsymbol{\alpha}_1$ $R_2^i$ $\boldsymbol{\alpha}_2 \rightarrow \boldsymbol{\alpha}_2$ $R_2$ $\boldsymbol{\alpha}_1$, because K ($\boldsymbol{\alpha}_1$) = K ($\boldsymbol{\alpha}_2$). Thus, there is always another $R_2^j$ such that $\boldsymbol{\alpha}_2$ $R_2^j$ $\boldsymbol{\alpha}_1$

- **Equivalence $\boldsymbol{\alpha_1}$ $\mathbf{R_2^{51}}$ $\boldsymbol{\alpha_2}$, $\boldsymbol{\alpha_1}$ = (α 0 0), $\boldsymbol{\alpha_2}$ = [(α+9)/2 (19-α)/2 5], α=5, 7, 9, w=$\dot{2}$**

  Indeed, under [18], $K_{31}(\boldsymbol{\alpha}) = (\alpha_1^1\text{-}1 \;\; 10\text{-}\alpha_1^1 \;\; \underbrace{0 \;\ldots\ldots}_{h-2} \;\; 0)$

  under [8a] with [8a.1], $K_1(\alpha_2^1 \; \alpha_2^2 \ldots \alpha_2^h) = (2\,\alpha_2^1\text{-}10 \;\; 2\,\alpha_2^2\text{-}9 \;\ldots\; 2\,\alpha_2^{h-1}\text{-}9 \;\; 2\,\alpha_2^h\text{-}10)$
  $K_{31}(\boldsymbol{\alpha}_1) = K_1(\boldsymbol{\alpha}_2) \rightarrow \alpha_1^1\text{-}1 = 2\,\alpha_2^1\text{-}10$; $10\text{-}\alpha_1^1 = 2\,\alpha_2^2\text{-}9$
  $2\,\alpha_2^s\text{-}9 = 0$, $s \geq 3$, ... h-1; $2\,\alpha_2^h\text{-}10 = 0$. These conditions are only possible if h = 3, yielding
  $\alpha_2^1 = (\alpha_1^1+9)/2$ , $\alpha_2^2 = (19\text{-}\alpha_1^1)/2$, $\alpha_2^3 = 5$, $\alpha_1^1 = 5, 7, 9$
  500 $R_2^{51}$ 775, 700 $R_2^{51}$ 865, 900 $R_2^{51}$ 955

- **Equivalence $\boldsymbol{\alpha_1}$ $\mathbf{R_2^{52}}$ $\boldsymbol{\alpha_2}$, $\boldsymbol{\alpha_1}$ = (α 0 0), $\boldsymbol{\alpha_2}$ = [(10-(α/2) 4+(α/2) 5], α = 2, 4, 6, w=$\dot{2}$**

  It results from using $K_{32}$ ($\boldsymbol{\alpha}_1$) instead of $K_{31}$ ($\boldsymbol{\alpha}_1$)
  And from this 200 $R_2^{52}$ 955, 400 $R_2^{52}$ 865, 600 $R_2^{52}$ 775
  Note that classes 100, 300 and 800 are excluded from $R_2^{52}$. Class 100 does not have an equivalent $R_2$ different from itself and 300 $R_2^3$ 800.
  Equivalences $R_2^{51}$ and $R_2^{52}$ do not operate on numbers with w-digits, w = $\dot{2}$+1.

- In $A_w$, w=$\dot{2}$+1, there are specific equivalences $R_2$, unshared with w=$\dot{2}$. But such $R_2$ vary with w, and thus general expressions do not exist. Table 5 shows all those existing in $A_7$ and $A_5$. Some of these equivalences, such as $e_2^{61}$, also generate equivalences $R_2$ from the same sequence, e.g. 977 $R_2$ 973. Furthermore, if the domains of existence of $K_i$ and $K_j$ have classes (α ß γ) in common, they generate equivalences $R_1$.

Table 5. Equivalences $R_2$ in $A_w$, w=$\dot{2}$+1, not found in w = $\dot{2}$

| All equivalences $R_2$ unique to $A_7$ (1) | | |
|---|---|---|
| Equivalences $R_2$ | Valid for | Derives from $K_i(\boldsymbol{\alpha}_1) = K_j(\boldsymbol{\alpha}_2)$ |
| $e_2^{61}$ (α ß γ) = [α (ß+γ)/2 10-ß] | 986 $R_2$ 972, 975 $R_2$ 963, 875 $R_2$ 863 | $K_1$, $K_7$ |
| $e_2^{62}$ (α ß γ) = [(α+ß+1)/2 ß+γ-α+2 11-α] | 986 $R_2$ 972, 985 $R_2$ 962, 966 $R_2$ 852, 875 $R_2$ 863 | $K_1$, $K_4$ |
| $e_2^{63}$ (α ß γ) = ( α 10-γ 11-ß) | 981 $R_2$ 993, 971 $R_2$ 994, 972 $R_2$ 984, 961 $R_2$ 995, 962 $R_2$ 985, 872 $R_2$ 884, 862 $R_2$ 885, 863 $R_2$ 875, 763 $R_2$ 775 | $K_4$, $K_{13}$ |
| $e_2^{64}$ (α ß γ) = (α 10-γ ß+2) | 931 $R_2$ 995, 921 $R_2$ 994, 911 $R_2$ 993, 932 $R_2$ 985, 922 $R_2$ 984, 832 $R_2$ 885, 822 $R_2$ 884, 833 $R_2$ 875, 733 $R_2$ 775 | $K_5$, $K_{13}$ |
| $e_2^{65}$ (α ß 0) = (11-ß α+ß-3 21-2α-ß) | 720 $R_2$ 965, 630 $R_2$ 866 | $K_{27}$, $K_1$ |
| $e_2^{66}$ (α ß 0) = (11-ß 9-α/2 13-α-ß) | 630 $R_2$ 864 | $K_{27}$, $K_7$ |
| $e_2^{67}$ (α ß 0) = (α+ß 18-2α-ß ß) | 610 $R_2$ 751, 520 $R_2$ 762 | $K_{27}$, $K_6$ |
| $e_2^{68}$ (α ß 0) = (20-2α-ß 10-α-ß ß) | 710 $R_2$ 521, 620 $R_2$ 622 | $K_{27}$, $K_{11}$ |
| $e_2^{69}$ (α ß 0) = (11-ß α+ß-3 12-α) | 730 $R_2$ 875, 720 $R_2$ 965, 630 $R_2$ 866 | $K_{26}$, $K_1$ |
| $e_2^{70}$ (α ß 0) = (α 19-α-ß α+2ß-10) | 930 $R_2$ 975, 840 $R_2$ 876 | $K_{25}$, $K_1$ |
| $e_2^{71}$ (α ß 0) = (ß+1 α-4 ß-α+3) | 880 $R_2$ 943, 770 $R_2$ 833 | $K_{23}$, $K_{12}$ |
| $e_2^{72}$ (α ß 0) = [(α+9)/2 20-2ß 10-ß] | 770 $R_2$ 863 | $K_{23}$, $K_4$ |
| $e_2^{73}$ (α ß 0) = (5+ß/2 21-α-ß 10-α) | 960 $R_2$ 862 | $K_{22}$, $K_4$ |
| All equivalences $R_2$ unique to $A_5$ (1) | | |
| $e_2^{5}$ (α ß) = [1/2 (α+ß+1) 11-α] | 85 $R_2$ 73, 96 $R_2$ 82 | $K_1$, $K_{17}$ (2) |
| $e_2^{6}$ (α 0) = [1/2 (11-α) α] | 30 $R_2$ 43 | $K_{26}$, $K_{14}$ |
| $e_2^{7}$ (α 0) = (11-α 2α-3) | 40 $R_2$ 75 | $K_{26}$, $K_1$ |
| (2) Notation for functions $K_i$ in $A_5$ is independent to that one in $A_7$. The permutations associated to $K_i$ in $A_7$ are shown in Table 2. The permutations associated to $K_i$ in $A_3$ are: $P_1$ (1234), $P_{17}$ (1342), $P_{14}$ (4312), $P_{26}$ (2 1) | | |

(1) There are no more different new $R_2$, except for those resulting from the transitivity (product) between $R_2$.

### 6.4.3 Equivalence classes $C_2$

Equivalence classes $C_2$, which include numbers whose second transformations match, are an essential element to understand the architecture of transformation trees. As an analogy, were parametric classes $C_1$ the leaves of a tree, classes $C_2$ would be the little branches that group the leaves together. The existence of equivalences $R_2$ contributes to new branches in the transformation trees by means of cycles.

In Table 6 equivalence classes $C_2$ in $A_6$ are shown. Some classes $C_1$ do not have equivalent $R_2$, such as 555 $\in C_2^{13}$, but others have up to 5 equivalent. Such is the case of $C_2^{74}$: any class 863, 833, 762, 732 and 332 belonging to it will become 643. And any number belonging to such classes will return the same image (642654) after 2 transformations. In this example 863 = $\boldsymbol{\alpha}_{c2}$, 643 = $\alpha_{c3}$ and 642654 = $n_{c4}$.

Graph 2 shows the transformations of classes $C_2$ in $A_6$. When compared to Graph 1, the structure of transformation trees stands out.

Table 6. Equivalence classes $C_2$ in $A_6$

| | Tree A. Branch A | | | | | | | | | | | |
|---|---|---|---|---|---|---|---|---|---|---|---|---|
| Equivalence classes | $C_2^1$ | $C_2^2$ | $C_2^3$ | $C_2^4$ | $C_2^5$ | $C_2^6$ | $C_2^7$ | $C_2^8$ | $C_2^9$ | $C_2^{10}$ | $C_2^{11}$ | $C_2^{12}$ |
| Classes $C_1$ = (α ß γ) included | 966<br>964<br>661<br>631<br>431 | 100 | 832<br>862 | 655<br>554<br>544 | 955<br>900<br>551<br>541<br>200 | 753<br>743 | 810<br>210 | 965<br>531 | 985<br>511 | 970<br>930 | 870<br>830 | 877<br>873<br>772<br>722<br>322 |
| Image of $C_2$ | 832 | 900 | 753 | 210 | 810 | 531 | 970 | 830 | 870 | 873 | 772 | 654 |
| $\boldsymbol{\alpha}_{ci}$ and D | 861;14 | 863;14 | 863;14 | 643;14 | 643;14 | 643;14 | 421;14 | 421;14 | 421;14 | 852;14 | 852;14 | 751;14 |

| $C_2^{13}$ | $C_2^{14}$ | $C_2^{15}$ | $C_2^{16}$ | $C_2^{17}$ | $C_2^{18}$ | $C_2^{19}$ | $C_2^{20}$ | $C_2^{21}$ | $C_2^{22}$ | $C_2^{23}$ | $C_2^{24}$ |
|---|---|---|---|---|---|---|---|---|---|---|---|
| 555 | 990<br>910<br>110 | 710<br>310 | 654<br>644<br>444 | 996<br>994 | 981<br>911<br>111 | 961<br>931 | 987<br>983<br>711<br>311 | 988<br>982<br>881<br>811<br>211 | 977<br>973<br>771<br>721<br>321 | 876<br>874<br>622<br>422 | 854<br>844<br>652<br>642<br>442 |
| 111 | 981 | 961 | 311 | 982 | 977 | 973 | 874 | 876 | 854 | 652 | 621 |
| 841;14 | 841;14 | 841;14 | 841;14 | 841;14 | 861;7 | 861;7 | 861;7 | 861;7 | 863;7 | 863;7 | 643;7 |

| Tree A. Branch B | | | | | | | | | | | | | |
|---|---|---|---|---|---|---|---|---|---|---|---|---|---|
| $C_2^{25}$ | $C_2^{26}$ | $C_2^{27}$ | $C_2^{28}$ | $C_2^{29}$ | $C_2^{30}$ | $C_2^{31}$ | $C_2^{32}$ | $C_2^{33}$ | $C_2^{34}$ | $C_2^{35}$ | $C_2^{36}$ | $C_2^{37}$ | $C_2^{38}$ |
| 995 | 980<br>920 | 888<br>882 | 770<br>730<br>330 | 666<br>664 | 775<br>600<br>500 | 963<br>933<br>761<br>731<br>331 | 540 | 665 | 800<br>300 | 853<br>843<br>752<br>742 | 660<br>640<br>440 | 720<br>320 | 954<br>944<br>651<br>641<br>441 |
| 980 | 882 | 775 | 761 | 331 | 540 | 843 | 640 | 320 | 720 | 641 | 651 | 860 | 821 |
| 751;14 | 841;14 | 861;14 | 863;14 | 863;14 | 863;14 | 643;14 | 643;14 | 643;14 | 643;14 | 421;14 | 421;14 | 421;14 | 852;14 |

## Table 6. (Continuation)

| $C_2^{39}$ | $C_2^{40}$ | $C_2^{41}$ | $C_2^{42}$ | $C_2^{43}$ | $C_2^{44}$ | $C_2^{45}$ | $C_2^{46}$ | $C_2^{47}$ | $C_2^{48}$ | $C_2^{49}$ | $C_2^{50}$ | $C_2^{51}$ |
|---|---|---|---|---|---|---|---|---|---|---|---|---|
| 880<br>820<br>220 | 860<br>840 | 951<br>941 | 972<br>922<br>871<br>821<br>221 | 763<br>733<br>333 | 971<br>921 | 865<br>700<br>532<br>400 | 765<br>533 | 975<br>521 | 630<br>430 | 850 | 750 | 754<br>744<br>653<br>643<br>443 |
| 871 | 763 | 971 | 865 | 533 | 975 | 630 | 430 | 850 | 750 | 754 | 653 | 421 |
| 852;14 | 852;14 | 751;7 | 751;7 | 751;7 | 841;7 | 841;7 | 841;7 | 861;7 | 861;7 | 863;7 | 863;7 | 643;7 |

| Tree A. Branch C | | | | | | | | | | | | | |
|---|---|---|---|---|---|---|---|---|---|---|---|---|---|
| Equivalence classes | $C_2^{52}$ | $C_2^{53}$ | $C_2^{54}$ | $C_2^{55}$ | $C_2^{56}$ | $C_2^{57}$ | $C_2^{58}$ | $C_2^{59}$ | $C_2^{60}$ | $C_2^{61}$ | $C_2^{62}$ | $C_2^{63}$ | $C_2^{64}$ |
| Classes $C_1 = (\alpha\ \beta\ \gamma)$ incluided | 999<br>991 | 998<br>992 | 997<br>993 | 875<br>522 | 986<br>984<br>611<br>411 | 776<br>774 | 650 | 950 | 777<br>773 | 872<br>822<br>222 | 542<br>552<br>855 | 755<br>553<br>543 | 610<br>410 |
| Image of $C_2$ | 997 | 986 | 984 | 650 | 872 | 542 | 552 | 855 | 553 | 755 | 610 | 410 | 952 |
| $\boldsymbol{\alpha}_{ci}$ and D | 861;7 | 863;7 | 863;7 | 643;7 | 643;7 | 421;7 | 421;7 | 421;7 | 421;7 | 421;7 | 852;7 | 852;7 | 751;7 |

| Tree A. Cycle zone | | | | | | | | | | |
|---|---|---|---|---|---|---|---|---|---|---|
| $C_2^{65}$ | $C_2^{66}$ | $C_2^{67}$ | $C_2^{68}$ | $C_2^{69}$ | $C_2^{70}$ | $C_2^{71}$ | $C_2^{72}$ | $C_2^{73}$ | $C_2^{74}$ | $C_2^{75}$ |
| 520 | 530 | 510 | 620 420 | 886 884 | 863<br>833<br>762<br>732<br>332 | 976 974<br>621 421 | 852 842 | 953 943<br>751 741 | 952 942<br>851 841 | 962 932<br>861 831 |
| 842 | 741 | 943 | 851 | 762 | 643 | 852 | 751 | 841 | 861 | 863 |
| 421;7 | 852;7 | 852;7 | 751;7 | 861;7 | 863;7 | 421;7 | 852;7 | 751;7 | 841;7 | 861;7 |

| Tree B | | | | | | Tree C |
|---|---|---|---|---|---|---|
| $C_2^1$ | $C_2^2$ | $C_2^3$ | $C_2^4$ | $C_2^5$ | $C_2^6$ | $C_2^1$ |
| 885 | 887 883 | 760 740 | 766 764<br>663 633<br>433 | 960 940 | 866 864<br>662 432 | 550 |
| 760 | 764 | 662 | 432 | 864 | 632 | 550 |
| 632;3 | 632;3 | 632;2 | 632;2 | 632;2 | 632;1 | 550;0 |

* D stands for the number of necessary transformations so that any class $C_1$ belonging to a given $C_2^i$ becomes the stated class

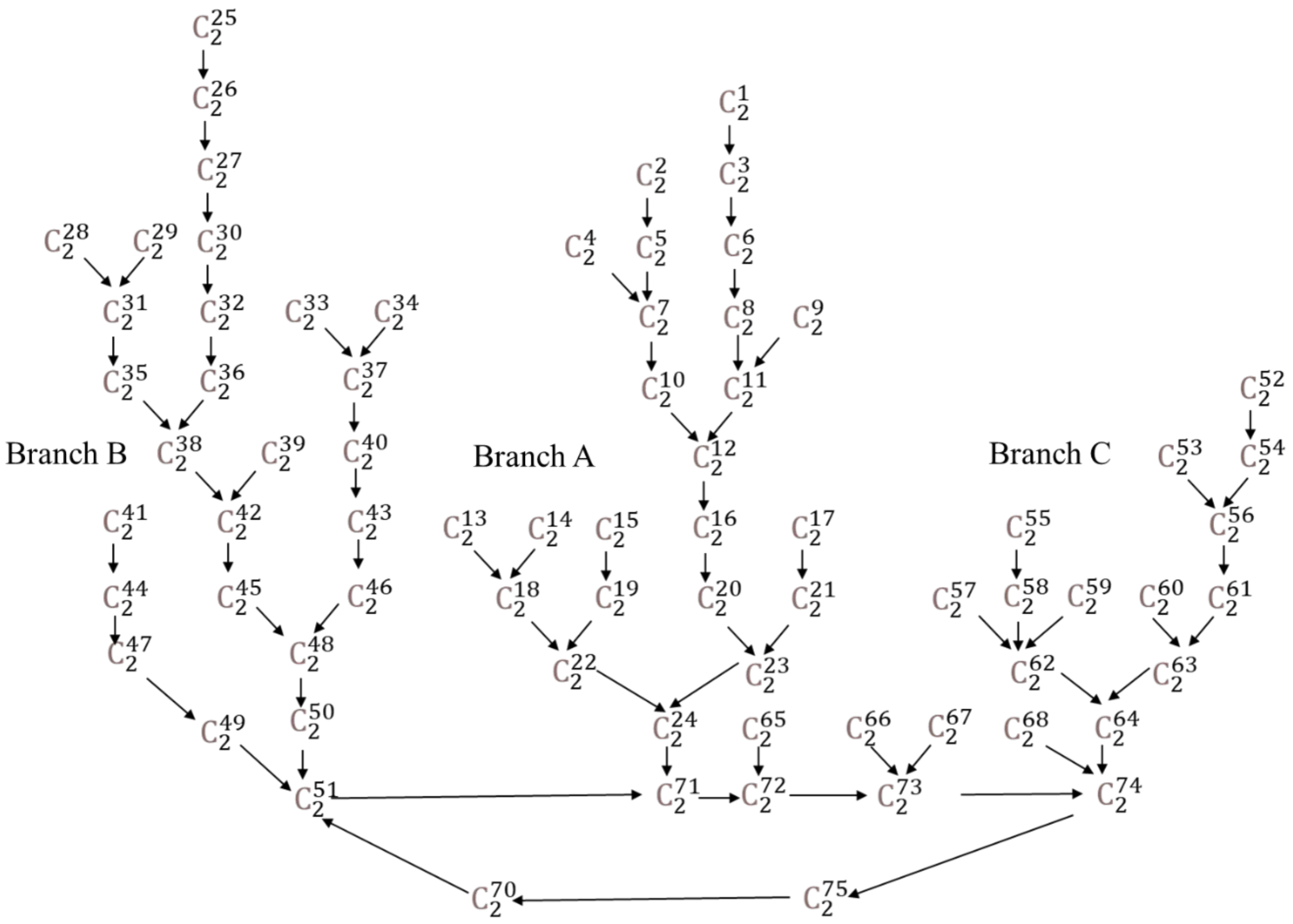

Tree A

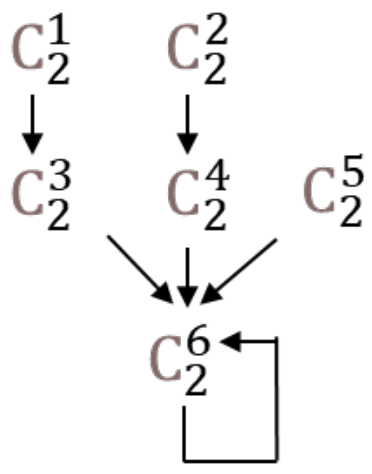

Tree B

Tree C

Graph 2. Transformation trees of equivalence classes $C_2$ in $A_6$

### 6.5 Groups of equivalence $R_2$

Transpositions in those sequences compatible with the transformation functions should in theory form subgroups of symmetric group $S_h$ of all h-element possible permutations. This should be true for the product of those equivalences in Set I based on [4a] or [4b], those in Set II based on [5a] or [5b] and those in Set III based on [6].

However, generally speaking it is not true for sets I and II. This is due to the fact that as h increases, some multiple transpositions are forbidden by [3]. Such is the case of $n_2^1$ and $n_2^2$ discussed in section 6.4.1. Specifically, for $A_6$ and $A_7$ the product of equivalences $e_2$ and transpositions $n_2$ of Set I are shown in Table 7. Apart from the general associativity in the products of transpositions, as well as the existence of the neutral element $e_2^0$ it is also worth noting the involutory $e_2^i \times e_2^i = e_2^0 \leftrightarrow (e_2^i) = (e_2^i)^{-1}$, and abelian $e_2^i \times e_2^j = e_2^j \times e_2^i$ features. But the product is no longer a closed binary relation between equivalences, since the product of some of them returns $n_2^6$ or $n_2^7$, which although being transpositions, are not equivalences.

As h decreases, the forbidden general multiple transpositions $n_2^1$ and $n_2^2$ disappear, and all transpositions generate equivalences. This happens when there are 2 non-null parameters –classes $(\alpha_1\ \alpha_2\ \overbrace{0 \ldots 0}^{h-2})$ from Set II and $(\alpha_1\ \alpha_2)$ from Set I with h = 2. Table 7 (Group II) shows the product of equivalences for $A_6$ and $A_7$. The same is true for the four equivalences of Set I in $A_4$ and $A_5$. The product of the equivalence in these cases has an isomorph group structure of Klein group. In addition, it is a subgroup of the symmetric group $S_4$ of all 4-element possible permutations.

Table 7. Product of equivalences $R_2$ in $A_6$ and $A_7$ based on the transposition of digits from the same sequence (f x g) (**α**) = f [g (**α**)]. (Equivalences $e_2$ are defined in Table 4).

| Group I | | | | | | | | |
|---|---|---|---|---|---|---|---|---|
| f\g | $e_2^0$ | $e_2^1$ | $e_2^2$ | $e_2^3$ | $e_2^4$ | $e_2^5$ | $n_2^6$ | $n_2^7$ |
| $e_2^0$ | $e_2^0$ | $e_2^1$ | $e_2^2$ | $e_2^3$ | $e_2^4$ | $e_2^5$ | $n_2^6$ | $n_2^7$ |
| $e_2^1$ | $e_2^1$ | $e_2^0$ | $e_2^4$ | $e_2^5$ | $e_2^2$ | $e_2^3$ | $n_2^7$ | $n_2^6$ |
| $e_2^2$ | $e_2^2$ | $e_2^4$ | $e_2^0$ | $n_2^6$ | $e_2^1$ | $n_2^7$ | $e_2^3$ | $e_2^5$ |
| $e_2^3$ | $e_2^3$ | $e_2^5$ | $n_2^6$ | $e_2^0$ | $n_2^7$ | $e_2^1$ | $e_2^2$ | $e_2^4$ |
| $e_2^4$ | $e_2^4$ | $e_2^2$ | $e_2^1$ | $n_2^7$ | $e_2^0$ | $n_2^6$ | $e_2^5$ | $e_2^3$ |
| $e_2^5$ | $e_2^5$ | $e_2^3$ | $n_2^7$ | $e_2^1$ | $n_2^6$ | $e_2^0$ | $e_2^4$ | $e_2^2$ |
| $n_2^6$ | $n_2^6$ | $n_2^7$ | $e_2^3$ | $e_2^2$ | $e_2^5$ | $e_2^4$ | $e_2^0$ | $e_2^1$ |
| $n_2^7$ | $n_2^7$ | $n_2^6$ | $e_2^5$ | $e_2^4$ | $e_2^3$ | $e_2^2$ | $e_2^1$ | $e_2^0$ |

| Group II | | | | |
|---|---|---|---|---|
| f\g | $e_2^0$ | $e_2^8$ | $e_2^9$ | $e_2^{10}$ |
| $e_2^0$ | $e_2^0$ | $e_2^8$ | $e_2^9$ | $e_2^{10}$ |
| $e_2^8$ | $e_2^8$ | $e_2^0$ | $e_2^{10}$ | $e_2^9$ |
| $e_2^9$ | $e_2^9$ | $e_2^{10}$ | $e_2^0$ | $e_2^8$ |
| $e_2^{10}$ | $e_2^{10}$ | $e_2^9$ | $e_2^8$ | $e_2^0$ |

| Group III | | |
|---|---|---|
| f\g | $e_2^0$ | $e_2^{11}$ |
| $e_2^0$ | $e_2^0$ | $e_2^{11}$ |
| $e_2^{11}$ | $e_2^{11}$ | $e_2^0$ |

When the number of non-null parameters decreases to 1, a cyclical subgroup isomorphic to $\mathbb{Z}_2$ appears. Such is the case of the general equivalence $e_2^3$ (section 6.4.1). Table 7 (Group III) shows the products for $A_6$ and $A_7$. This situation occurs as well in $A_2$ and $A_3$.

Naturally, any involutory equivalence, along with the identical, form a group isomorphic to $\mathbb{Z}_2$. Such is the case of $e_2^{11}$ in $A_6$, which in this case does not exist in $A_7$. However, its derivative functions (Table 4), along with the equivalence itself, do not have a subgroup structure, since the product of the former generates equivalences of Set I.

### 6.6 Symmetries $R_3$ and higher

From [28b] and [21] it results

$\boldsymbol{\alpha}_1 R_r \boldsymbol{\alpha}_2 \leftrightarrow \boldsymbol{\alpha}'_1 R_{r-1} \boldsymbol{\alpha}'_1$, $\boldsymbol{\alpha}'_1 = K_i(\boldsymbol{\alpha}_1)$, $\boldsymbol{\alpha}'_2 = K_j(\boldsymbol{\alpha}_2)$ [33]

Specifically: $\boldsymbol{\alpha}_1 R_3 \boldsymbol{\alpha}_2 \leftrightarrow \boldsymbol{\alpha}'_1 R_2 \boldsymbol{\alpha}'_2$ [33.1]

One of the approaches to obtain the algebraic expressions of equivalences $R_3$ is based on this expression. If an equivalence $R_2$ is imposed on the images of two classes $\boldsymbol{\alpha}_1$ and $\boldsymbol{\alpha}_2$, $K_i(\boldsymbol{\alpha}_1)$ and $K_j(\boldsymbol{\alpha}_2)$, the result will be an algebraic relation between $\boldsymbol{\alpha}_1$ and $\boldsymbol{\alpha}_2$ so that

$\boldsymbol{\alpha}_1 R_3 \boldsymbol{\alpha}_2$: $(e_2 \times K_i)(\boldsymbol{\alpha}_1) = K_j(\boldsymbol{\alpha}_2)$ [34]

As a few examples we introduce the following equivalences

**Equivalence $\boldsymbol{\alpha}_1$ $R_3^1$ $\boldsymbol{\alpha}_2$, $\boldsymbol{\alpha}_1 = (\alpha^1\ \alpha^2 ... \alpha^h)$, $\boldsymbol{\alpha}_2 = (15\text{-}\alpha^1\ \alpha^2 ... \alpha^h)$**

$w = \dot{2}$; $6 \leq \alpha^1 \leq 9$; $\alpha^1 \geq \alpha^2+1$, $\alpha^1 + \alpha^2 \leq 14$; $5 \leq \alpha^s \leq 6$, $s \geq 2$;

For example, in $A_6$

$(\alpha\ \beta\ \gamma)$ $R_3^1$ $(15\text{-}\alpha\ \ \beta\ \ \gamma)$

Which provides the following equivalences

Tree A: 955 $R_3^1$ 655, 865 $R_3^1$ 765, 855 $R_3^1$ 755

Tree B: 866 $R_3^1$ 766

or in terms of classes $C_2$, for [29]

$C_2^5$ $R_3$ $C_2^4$, $C_2^{45}$ $R_3$ $C_2^{46}$, $C_2^{65}$ $R_3$ $C_2^{66}$, $C_2^6$ $R_3$ $C_2^4$ (Tree B)

**Equivalence $\boldsymbol{\alpha}_1$** $R_3^2$ **$\boldsymbol{\alpha}_2$**, $\boldsymbol{\alpha}_1 = (\alpha^1\ \alpha^2 \ldots \alpha^h)$, $\boldsymbol{\alpha}_2 = (\alpha^1\ \alpha^2 \ldots \alpha^{h-1}\ \ 15\text{-}\alpha^h)$

$w = \dot{2}$, $6 \leq \alpha^h \leq 9$, $\alpha^{h-1}+\alpha^h \geq 15$, $\alpha^1 \geq \alpha^2+1$

Derives from $e_2^{13} \times K_1(\boldsymbol{\alpha}_1) = K_1(\boldsymbol{\alpha}_2)$

For example, in $A_6$

$(\alpha\ \beta\ \gamma)$ $R_3^2$ $(\alpha\ \ \beta\ \ 15\text{-}\gamma)$ valid for

988 $R_3^2$ 987, 987 $R_3^2$ 988 $\rightarrow$ $C_2^{21}$ $R_3$ $C_2^{20}$

**Equivalence $\boldsymbol{\alpha}_1$** $R_3^3$ **$\boldsymbol{\alpha}_2$**, $\boldsymbol{\alpha}_1 = (\alpha^1\ \alpha^2 \ldots \alpha^h)$, $\boldsymbol{\alpha}_2 = (\alpha^1\ \alpha^2 \ldots \alpha^{h-1}\ \ 28\text{-}2\alpha^{h-1}\text{-}\alpha^h)$

$w = \dot{2} + 1$, $5 < \alpha^1 \geq \alpha^2+1$, $5 \leq \alpha^s \leq 8$, $s>1$; $\alpha^s + \alpha^{s+1} \geq 9$, $s=1,2 \ldots h\text{-}1$, $3\ \alpha^{h-1} + \alpha^h \geq 28$

$2\alpha^{h-1} + \alpha^h \leq 23$

Derives from $e_2^{13} \times K_1(\boldsymbol{\alpha}_1) = K_1(\boldsymbol{\alpha}_2)$

For example, in $A_7$

$(\alpha\ \beta\ \gamma)$ $R_3^3(\alpha\ \ \beta\ \ 28\text{-}2\beta\text{-}\gamma)$ valid for the new equivalences

987 $R_3^3$ 985, 985 $R_3^3$ 987 and for $R_1$ 986, 977 y 877

**Equivalence $\boldsymbol{\alpha}_1$** $R_3^4$ **$\boldsymbol{\alpha}_2$**, $\boldsymbol{\alpha}_1 = (\alpha\ 0\ 0)$, $\boldsymbol{\alpha}_2 = [1/2\ (21\text{-}\alpha)\ \ 5\ \ 5]$

Equivalence $R_3$ is often limited to a couple of classes. Such is the case of

$e_2^8 \times K_{31}(\boldsymbol{\alpha}_1) = K_7(\boldsymbol{\alpha}_2)$ in $A_6$.

$e_2^8 \times K_{31}(\alpha_1\ 0\ 0) = (11\text{-}\alpha_1\ \ 10\text{-}\alpha_1\ \ 0) = (2\alpha_2\text{-}10\ \ \beta_2\text{-}\gamma_2\text{+}1\ \ \beta_2\text{-}\gamma_2) = K_7(\boldsymbol{\alpha}_2)$

which entails $\gamma_2=\beta_2$, $\alpha_1=9$, $\alpha_2=6$ and since $K_7$ demands $5 \leq \beta \leq 8$ y $2 \leq \gamma \leq 5 \rightarrow \gamma = \beta = 5$

the equivalence is reduced to 900 $R_3^4$ 655

This situation becomes generalized when increasing the order of the equivalence. Thus, in $A_7$

**Equivalence** $\boldsymbol{\alpha}_1$ $R^1_4$ $\boldsymbol{\alpha}_2$, $\boldsymbol{\alpha}_1 = (\alpha\ ß\ \gamma)$, $\boldsymbol{\alpha}_2 = (\alpha \quad 48-4\alpha-ß+2\gamma \ \ \gamma)$

Derives from $e^3_3$ x $K_4$ ($\boldsymbol{\alpha}_1$) = $K_4$ ($\boldsymbol{\alpha}_2$), only valid for 981 $R^1_4$ 961
which answers another of the questions in the introduction:
5068069 $R_4$ 3071934 because p(5068069) = 981 and p(3071934) = 961

**Equivalence** 533 $R^1_7$ 621,

Indeed $K^6$ (533) = $K^6$ (621) = 864. Or equivalently, m=4687437, p(m)=533, n=4693554, p(n)=(621), $K^7$(m)=$K^7$(n)=8639532 →m $R_7$ n

However, the algebraic relation between 533 and 621, or between m and n, as in the two previous examples, is not extendible to other classes $\boldsymbol{\alpha}_1$ $R_7$ $\boldsymbol{\alpha}_2$.

As order r of equivalence $R_r$ increases, it is less interesting to establish the determining algebraic relation.

### 6.7 Symmetries and cyles

The transformation trees are articulated on the cycles, whether they have several links r or only one (transformation constants).

Equivalence classes $C_r$ group numbers with the same image $C_{r-1}$. In general, $\boldsymbol{\alpha}_1$ $R_s$ $\boldsymbol{\alpha}_2$ ↔ $\boldsymbol{\alpha}'_1$ $R_{s-1}$ $\boldsymbol{\alpha}'_2$. For example, classes $C_1$ 963, 933, 761, 731, and 331 all belong to $C^{31}_2$ because their image is the same, 843 (Table 6). Analogously, $C^4_4$ and $C^5_4$ belong to the same class $C_5$ because they have the same image $C^6_4$ (Graph 3). This relation forces the transformation tree to swallow the leaves and little branches as s increases, and therefore only the big branches and the trunk remain visible.

The absorption process occurs in the whole tree as well as in the cycles. The latter act as black holes that swallow those classes converging with them. Each class $C_s$ (ci) from the cycle groups with classes $C_s$ which have as an image $C_s$ ($c_{i+1}$) from the cycle, resulting in a new class $C_{s+1}$ (cj). The process continues until $C_s = C_{s+1}$ for all classes $C_s$ in the tree. This way, we reach the top of new equivalences $R_s = R_u$, where there are only as many classes $C_u$ as there are links r in the cycle, $C^i_u$, i = 1, 2,…r.

Each class_ $C^i_u$ belongs to the cycle and groups all classes $C_s$, s = 0, 1, 2, u-1, whose distance (number of transformations) to $\boldsymbol{\alpha}_{ci}$ is a multiple of r. Thus, the number of links in the cycle establishes the distance pattern for the subsequent groups.

For instance, Table 8 shows the 7 classes $C^i_u$ resulting from u = 13 transformations in tree A from $A_6$. As stated several times, tree A has a cycle of r = 7 links.
Class 661 ∈ $C^1_2 \subset C^1_4 \subset C^1_{13}$ and $\boldsymbol{\alpha}_{C1}$ = 861 ∈ $C^{75}_2 \subset C^{28}_4 \subset C^1_{13}$ (Table 9). 661 needs 14 (2 x r) transformations to become 861 –the same number of transformations that $C^1_2$ needs to reach $C^{75}_2$ (Graph 2) or $C^1_4$ to reach $C^{28}_4$ (Graph 3). Class 430 ∈ $C^{48}_2$ only requires r = 7 transformations

Table 8. Inclusion relations $C_2 \subset C_4 \subset C_{13}$ in $A_6$

| $C_{13}^1$ | | $C_{13}^2$ | | $C_{13}^3$ | | $C_{13}^4$ | | $C_{13}^5$ | | $C_{13}^6$ | | $C_{13}^7$ | |
|---|---|---|---|---|---|---|---|---|---|---|---|---|---|
| $C_2$ | $C_4$ | $C_2$ | $C_4$ | $C_2$ | $C_4$ | $C_2$ | $C_4$ | $C_2$ | $C_4$ | $C_2$ | $C_4$ | $C_2$ | $C_4$ |
| $C_2^1$ | $C_4^1$ | $C_2^2$ | $C_4^2$ | $C_2^4, C_2^5$ | $C_4^4$ | $C_2^7, C_2^8, C_2^9$ | $C_4^6$ | $C_2^{10}, C_2^{11}$ | $C_4^7$ | $C_2^{12}$ | $C_4^8$ | $C_2^{13}, C_2^{14}, C_2^{15}$ | $C_4^9$ |
| $C_2^{18}, C_2^{19}, C_2^{20}, C_2^{21}$ | $C_4^{11}$ | $C_2^3$ | $C_4^3$ | $C_2^6$ | $C_4^5$ | $C_2^{35}, C_2^{36}$ | $C_4^{20}$ | $C_2^{38}, C_2^{39}$ | $C_4^{22}$ | $C_2^{25}$ | $C_4^{13}$ | $C_2^{16}, C_2^{17}$ | $C_4^{10}$ |
| $C_2^{27}$ | $C_4^{15}$ | $C_2^{28}, C_2^{29}$ | $C_4^{16}$ | $C_2^{24}, C_2^{51}$ | $C_4^{37}$ | $C_2^{37}$ | $C_4^{21}$ | $C_2^{40}$ | $C_4^{23}$ | $C_2^{41}$ | $C_4^{24}$ | $C_2^{26}$ | $C_4^{14}$ |
| $C_2^{52}$ | $C_4^{29}$ | $C_2^{30}$ | $C_4^{17}$ | $C_2^{31}, C_2^{32}$ | $C_4^{18}$ | $C_2^{57}, C_2^{58}, C_2^{59}, C_2^{60}, C_2^{61}$ | $C_4^{33}$ | $C_2^{62}, C_2^{63}, C_2^{66}, C_2^{67}, \mathbf{C_2^{72}}$ | $\mathbf{C_4^{34}}$ | $C_2^{42}, C_2^{43}$ | $C_4^{25}$ | $C_2^{44}$ | $C_4^{26}$ |
| $C_2^{47}, C_2^{48}, C_2^{69}, \mathbf{C_2^{75}}$ | $\mathbf{C_4^{28}}$ | $C_2^{53}, C_2^{54}$ | $C_4^{20}$ | $C_2^{33}, C_2^{34}$ | $C_4^{19}$ | $C_2^{65}, \mathbf{C_2^{71}}$ | $\mathbf{C_4^{38}}$ | | | $C_2^{64}, C_2^{68}, \mathbf{C_2^{73}}$ | $\mathbf{C_4^{35}}$ | $C_2^{45}, C_2^{46}$ | $C_4^{27}$ |
| | | $C_2^{22}, C_2^{23}, C_2^{49}, C_2^{50}, \mathbf{C_2^{70}}$ | $\mathbf{C_4^{12}}$ | $C_2^{55}$ | $C_4^{31}$ | | | | | | | $\mathbf{C_2^{74}}$ | $\mathbf{C_4^{36}}$ |
| | | | | $C_2^{56}$ | $C_4^{32}$ | | | | | | | | |
| *Class $C_{13}^i$ includes class $\boldsymbol{\alpha}_{ci}$. $C_{13+r}^i = C_{13}^i$, r>0. Classes $C_2$ and $C_4$ containing $\boldsymbol{\alpha}_{ci}$ are marked in bold. | | | | | | | | | | | | | |

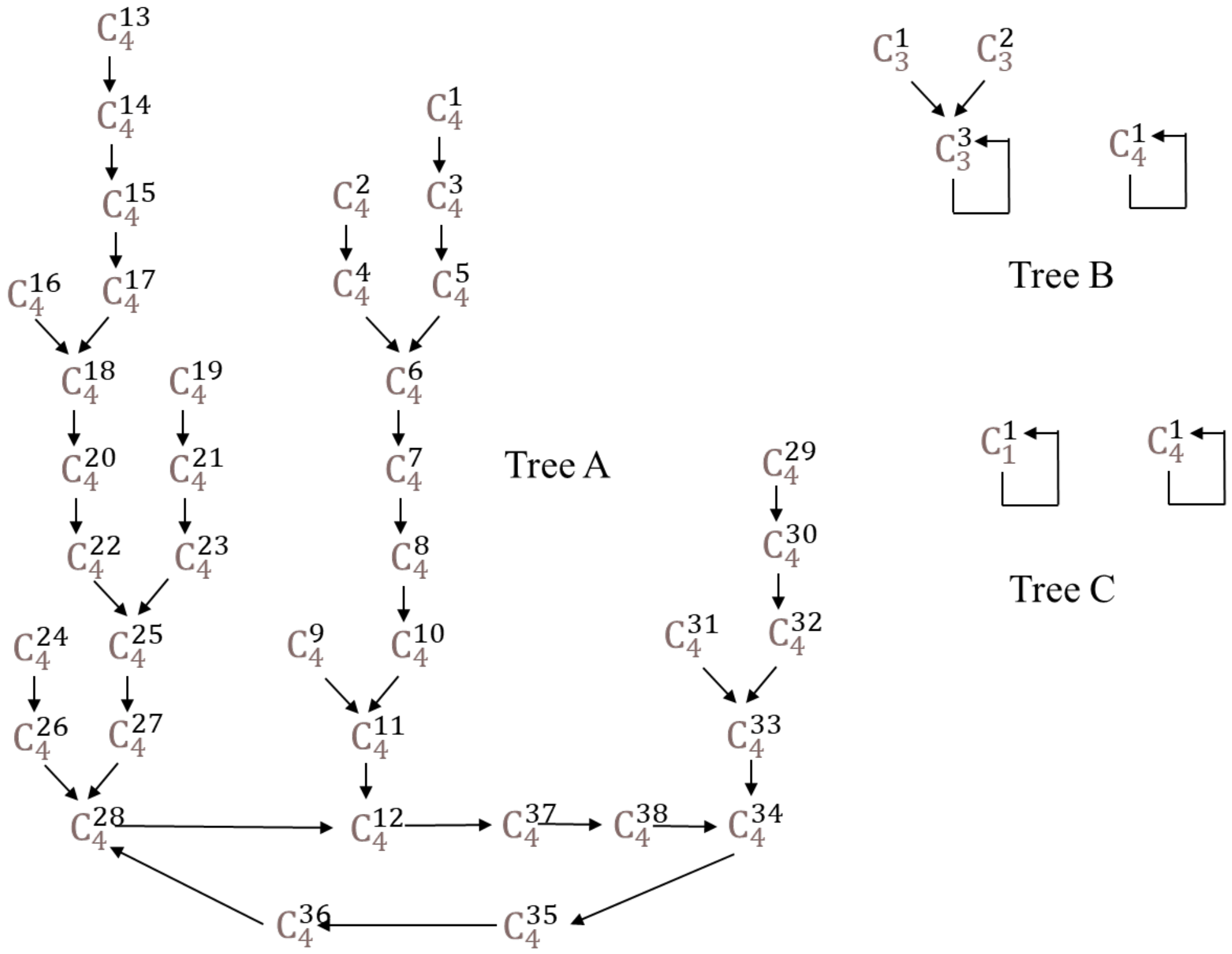

Graph 3. Transformations between equivalence classes $C_4$ in $A_6$. ($C_3$ are included in tree B, while $C_1$ appear in tree C)

## 7 Discussion

We understand symmetry as invariance during transformation. This has been the aim of these papers –trying to understand why sets of different numbers return identical images. Not only with a single transformation, but with any number r of them. To this end, binary equivalence relations of order r have been used. Two numbers m and n present an equivalence relation m $R_r$ n if and only if $K^r(m) = K^r(n)$. This simple relationship allows the numbers to be grouped into $C_r$ equivalence classes such that all their numbers give the same transform after r transformations. The use of these equivalences together with the complete parameterization of the transformation by means of $K_i$ functions has enabled us to find algebraic relations which explain the invariance of transformations.

The first step was to parameterize numbers. Any number n can have some parameters $\boldsymbol{\alpha}$ which determine its image. We established the general functions f that conduct such transformation $f(\boldsymbol{\alpha}) = n'$, n' being the image of n under Kaprekar's routine. These functions are similar to those established by Prichett et al. (1981). From here, equivalence classes of order 1 arise –those consisting of numbers with equal parameters and give the same transformed. The transformed numbers are multiples of 9 that satisfy some strict requirements. We have analyzed this aspect in detail, as it conditions the existence of symmetries.

The second step was to determine the parameters $\boldsymbol{\alpha}$ of the transformed number. Here lies one of the main problems with transformation. The transformed number need not have its digits sorted in ascending or descending order, which is a prerequisite to determine its parameters. The permutation that arranges the digits must therefore be established. However, this is no easy task, since the number of permutations grows factorially as the number of digits increases. Once both the permutation and the domain of existence have been established, **α**' is automatically defined. We have established the procedure to construct the functions $K_i$ (**α**) = **α** ', each one associated with a permutation. It is a novel contribution.

With these tools we were able to approach equivalences $R_r$ of a higher order. For this, the simple relation $\boldsymbol{\alpha}_1$ $R_r$ $\boldsymbol{\alpha}_2$ ↔ $\boldsymbol{\alpha}'_1$ $R_{r-1}$ $\boldsymbol{\alpha}'_2$ can be very useful. Starting from equivalences $R_1$, we can then move onto $R_2$, $R_3$, etc.

We showed that there are some $R_2$ –relations which must exist between numbers for their second transformation to match– that are universal. Those in sets I, II and III, valid for numbers with both an even number of digits (w = $\dot{2}$) and an odd number (w = $\dot{2}$+1). Others are only valid for w = $\dot{2}$. Thus, for example, the numbers whose parameters are $\boldsymbol{\alpha}_1 = (\alpha^1\ \alpha^2\ \ldots\alpha^h)$ give the same numbers after two transformations as those with parameters $\boldsymbol{\alpha}_2 = (10\text{-}\alpha^h\ \ 9\text{-}\alpha^{h-1}\ \ 9\text{-}\alpha^{h-2}\ldots 9\text{-}\alpha^3\ \ 9\text{-}\alpha^2\ \ 10\text{-}\alpha^1)$, $a^1 \geq a^2+1$. It is a $\boldsymbol{\alpha}_1$ $R_2$ $\boldsymbol{\alpha}_2$ valid for w= 2h, that is, for numbers with an even number of digits. It is an involutional equivalence.

The products of equivalences I and II have group algebraic structures only for low values of w. The equivalences of set I, for w = 4 or 5, while those of set II, for w = 6 or 7. Therefore, the group is isomorphic of Klein group. The equivalences of set III is universal, thus generating a group isomorphic to $\mathbb{Z}_2$. The reason why there are not any more groups lies in relation [3], which forces the parameters to be arranged. Many permutations violate this arrangement, resulting in transpositions without domain of existence in the product of equivalences.

The new higher equivalence classes demand increasingly narrow domains of existence. Their algebraic relation is often only valid for a few classes $C_i$. With this, they are less and less interesting as r increases.

The development of functions $K_i$ allowed us to approach another feature of transformations –the relation between cycles and symmetries. Transformation trees are articulated on cycles. They can be leafy or not. Thus, in $A_6$, tree A contains 201 classes $C_1$ and it is articulated on a 7-link cycle. Tree B contains 17 classes and it is articulated on the constant 631764, whose parameters are $\boldsymbol{\alpha}_E$ = 632. Tree C only contains class $\boldsymbol{\alpha}_E$ = 550, whose numeric constant is 549945. In $A_7$ the 219 existing parametric classes are grouped on an 8-link cycle. How many classes make up a tree will depend on the structure of functions $K_1$, or equivalently, on the compatibility between successive permutations. This compatibility varies with $A_w$.

Also dependent on each $A_w$ is the coexistence of cycles and constants –in $A_2$, there is a single cycle, in $A_3$ and $A_4$, a tree with a constant (a Kaprekar constant), in $A_5$, three cycles, in $A_6$, one cycle and two constants, in $A_7$, a single cycle, etc.

The structuration of the transformation trees in equivalence classes follows a grouping pattern mediated by the number of links in the cycle. This number establishes the "meter" of distances for the successive groupings of the equivalence classes.

In short, Kaprekar's routine raised the issue of the uniqueness of its constants. Prichett et al. (1981) showed that there are only two Kaprekar constants in base 10. Also, in this paper we have seen that there are several unique facts. The basic transformation applied to numbers in base 10 forces not only the transformed number to be a multiple of 9, but also its digits to satisfy some strict requirements ([4a.1 ], [4b.1]…). Such numeric and parametric restrictions ([3]) lead in many cases to uniqueness. The algebraization of Kaprekar's routine has proven as well the existence of universal relations. Such is the case of some $R_2$ symmetries and the existence of universal groups, such as the Klein group.

We hope that our contribution will help arouse interest in this and other numeric transformations.

## 8 References

- Kaprekar, D. 1949. Another solitaire game. Scripta Mathematica 15: 244-245.
- Gardner M. 1975. Mathematical games. Scientific American 232 (3): 112-117.
- Lapenta JF, Ludington AL, Prichett GD. 1979. An algorithm to determine self-producing r-digit g-adic integers. Journal of Mathematik. Band 310: 14.
- Prichett GD, Ludington AL, Lapenta JF. 1981. The determination of all decadic Kaprekar constants. Fibonacci Quarterly 19 (1): 45-52.
- Walden BL. 2005. Searching for Kaprekar's constants: algorithms and results. International Journal of Mathematics and Mathematical Sciences 2005 (18): 2999-3004.
- Nishiyama, Y. 2006. Mysterious number 6147. Plus.maths.org/content/mysterious-number-6147.
- Dolan S. 2011. A classification of Kaprekar constants. The Mathematical Gazette 95(534): 437-443.
- Hanover D. 2017. The base dependent behavior of Kaprekar's routine: A theoretical and computational study revealing new regularities. arXiv: 1710.06308v1 [math.GM].
- Yamagami A. 2018. On 2-adic Kaprekar constants and 2-digit Kaprekar distances. Journal of Number Theory 185: 257-280.
- Devlin P, Zeng T. 2020. Maximum distances in the four-digit Kaprekar process. arXiv: 2010.11756v1 [math.NT].
- Wang Z, Lu W. 2021. On 2-digit and 3-digit Kaprekar's routine. arXiv: 2101.09708 [math.NT].

**Acknowledgements**

To Dr. María José Díez, professor at the UPV, for her selfless help in the correction of this manuscript.